\documentclass[preprint]{article}
\usepackage{standalone}

\usepackage[utf8]{inputenc} 
\usepackage[T1]{fontenc}    
\usepackage{hyperref}       
\usepackage{url}            
\usepackage{booktabs}       
\usepackage{amsfonts}       
\usepackage{nicefrac}       
\usepackage{microtype}      
\usepackage{amsmath}
\usepackage{microtype}
\usepackage{graphicx}
\usepackage{subcaption}
\usepackage{booktabs} 
\usepackage{float}
\usepackage[algo2e,ruled,vlined,linesnumbered]{algorithm2e}
\usepackage{comment}
\usepackage{xcolor}
\usepackage{import}
\newtheorem{theorem}{Theorem}[section]

\newcommand{\mbb}{\mathbb}
\newcommand{\mbf}{\mathbf}
\newcommand{\mcl}{\mathcal}

\newcommand{\f}{\frac}

\newcommand{\T}{\textnormal}
\newcommand{\x}{\mathbf{x}}

\usepackage[normalem]{ulem} 
 
\title{Recursive Two-Step Lookahead Expected Payoff for Time-Dependent Bayesian Optimization}

\author{%
  S. Ashwin Renganathan, Jeffrey Larson and Stefan Wild\\
  Mathematics \& Computer Sciences\\
  Argonne National Laboratory\\
  Lemont, IL 60439 \\
  \texttt{srenganathan@anl.gov} \\
}

\begin{document}

\maketitle

\begin{abstract}
We propose a novel Bayesian method to solve the maximization of a
time-dependent expensive-to-evaluate oracle. We are interested in the decision
that maximizes the oracle at a finite time horizon, when relatively few noisy
evaluations can be performed before the horizon. Our recursive, two-step
lookahead expected payoff ($\texttt{r2LEY}$) acquisition function makes
nonmyopic decisions at every stage by maximizing the estimated expected value
of the oracle at the horizon. $\texttt{r2LEY}$ circumvents the evaluation of the
expensive multistep (more than two steps) lookahead acquisition function by
recursively optimizing a two-step lookahead acquisition function at every
stage; unbiased estimators of this latter function and its gradient are
utilized for efficient optimization. $\texttt{r2LEY}$ is shown to exhibit natural
exploration properties far from the time horizon, enabling accurate
emulation of the oracle, which is exploited in the final decision made at the
horizon. To demonstrate the utility of $\texttt{r2LEY}$, we compare it with
time-dependent extensions of popular myopic acquisition functions via both
synthetic and real-world datasets.
\end{abstract}

\section{Introduction}
\label{s:Introduction}
We consider the maximization of an expensive-to-evaluate oracle $f$ with
a finite budget of evaluations $q$. The inputs to $f$ are $\x$, the \emph{action}
parameters from a compact domain $\mcl{X}\subset \mbb{R}^d$, and $t$, the
\emph{context} from a (possibly infinite) set of contexts $\mcl{T}$.
The observations $y \in \mbb{R}$ made by observing $f$ at the action-context
pair $(\x, t)$ are assumed to be corrupted by additive stochastic noise $\epsilon$; that
is, $y = f(\x,t) + \epsilon$. Furthermore, we assume a 
context space
$\mcl{T}$ whose members have a unique ordering (e.g., time), and we are
interested in determining an optimum action at a given future context $T \in \mcl{T}$ using
only $q$ evaluations of $f$. Such a problem fundamentally differs
from general context-dependent $k$-armed bandit problems 
in that we are interested
solely in the action that maximizes payoff at context 
$t=T$ as opposed to maximizing the \emph{cumulative payoff} until $T$, as done
by, for example, \cite{krause2011contextual} and \cite{srinivas2009gaussian}.
Additionally, we assume that only one observation can be made per context (hereafter
``time'') $t$ and that the schedule of observations $\lbrace t_i \rbrace \forall~ i=1,\ldots,q$ is given.
The time-dependent maximization that we address arises, for example, in quantum
computing applications~\cite{zhu2018training}
where, in order
to find optimal tuning parameters $\x$ at time $T$, quantum circuit parameters
are tuned by using real-time noisy observations from the quantum device
$f(\x,t)$. 

The main challenges of our problem are the unknown structure of $f$
and the high costs associated with each evaluation of $f(\x,t)$.
Bayesian optimization (BO)~\cite{brochu2010tutorial, shahriari2015taking} with
Gaussian process (GP) priors~\cite{rasmussen:williams:2006}  suits the
specific
challenges posed by our problem,
where observations at each round $i$, 
$y_i =
f(\x_i, t_i) + \epsilon_i$, are made judiciously by leveraging information
gained from previous observations, as a means of coping with the high costs
of
each observation. The key idea 
is to specify GP prior distributions on the
oracle and the noise. That is, $f(\x,t) \sim \mcl{GP}\left(0, k( (\x, t), (\x',
t'))\right)$ and $\epsilon \sim \mcl{GP}(0,\sigma^2_\epsilon)$, where $k$ is a
covariance function, $\sigma^2_\epsilon$ is a constant noise variance and we
assume that $\mcl{X}$ is a hypercube with lower and upper bounds $lb$ and $ub$,
respectively. Note
that the covariance function captures the correlation between the observations
in the $(\x,t)$ space \emph{jointly}; here we use the product composite form
given by $k((\x, t), (\x', t')) = k_{\x}(\x,\x'; \theta_\x) \times k_t(t,t';
\theta_t) $, where $\theta_\x$ and $\theta_t$ 
parametrize the covariance functions for $\x$ and $t$, respectively, and 
$\Omega= \lbrace \theta_x, \theta_t, \sigma^2_\epsilon \rbrace$ are the unknown
GP hyperparameters that are estimated from data by maximizing the marginal
likelihood. The
posterior predictive distribution of the output $Y$
conditioned on available observations from the oracle, namely, $Y(\x, t |
\mcl{D}_n, \Omega) \sim \mcl{GP}(\mu_n(\x, t), \sigma^2_n (\x, t)),~ \mcl{D}_n = \lbrace
(\x_i, t_i), y_i \rbrace _{i=1}^n$,  
is then used as a surrogate model for $f$. Note that $\mu_n$ and $\sigma_n^2$
are the posterior mean and variance of the GP, respectively, where the
subscript $n$ implies the conditioning based on $n$ past observations. BO then
proceeds by defining an \emph{acquisition function} in terms of the GP posterior that is
optimized in lieu of the expensive $f$ to select the next
point $\x_{n+1}$,
and the process
continues recursively until either the budget $q$ is reached or the global
maximum of the
payoff function is realized; see Algorithm~\ref{a:BO}.

The standard acquisition functions in BO take a \emph{greedy} or \emph{myopic}
approach, where each decision in the sequence is made as though it were the
last, without accounting for the potential impact on the future decisions.
Several greedy acquisition functions have been proposed,
including the
probability of improvement (PI)~\cite{kushner1964new}, the expected improvement
(EI) \cite{jones1998efficient, mockus1978application}, and the GP upper
confidence bound (UCB)~\cite{srinivas2009gaussian}. Whereas in PI and EI the
acquisition function is a probabilistic measure of \emph{improvement}
over a user-specified target, in UCB it is an
optimistic estimate of the payoff.
From a finite-budget BO perspective,
such acquisition functions can be suboptimal~\cite{gonzalez2016glasses} and
moreover, optimizing them is guaranteed to attain the global optimum (under suitable
regularity conditions) only in the limit $q \rightarrow \infty$, e.g., see
~\cite{bull2011convergence, vazquez2010convergence}. For finite-budget
time-dependent
optimization problems with a target time horizon $T$, an acquisition that seeks
to optimize the \emph{longsighted} decision at $T$ is more appropriate.
In other words, we want an approach that makes decisions at each $t$ aimed at maximizing the
payoff at $T$. Such finite-budget BO strategies are sometimes referred as
\emph{lookahead} approaches, since the current decision is made by looking ahead
at future decisions.

In foundational work on finite-budget BO, Osborne \cite{osborne2010bayesian} showed that proper Bayesian reasoning
can be used to define an acquisition function where the $n+1$th
observation is made by marginalizing a loss function at the final ($q$th)
observation,  
over all the remaining $(\x_i, y_i)$, $\forall i=n+2,\ldots, q$
observations. However, the approach in \cite{osborne2010bayesian} is restricted
to a context-free setting and, furthermore, entailed defining the loss in terms
of the best observed value, which, as we will show, does not apply to
time-dependent problems.
Another way to solve the finite-budget BO problem is to sequentially choose points that are most \emph{informative}~\cite{cover2012elements} about the 
global
maximum/maximizer (or minimum/minimizer), e.g., see \cite{hennig2012entropy,
hernandez2014predictive, wang2017max}. However, such information-theoretic
approaches are known to involve an intractable form of the entropies, which
often need an approximation and 
are tailored
to seek a global optimum in a context-free setting (unlike our problem), and
hence their suitability for time-dependent BO is unknown.
Others have proposed
\emph{lookahead EI} in the context of finite-budget BO: Ginsbourger and Le Riche \cite{ginsbourger2010towards}
showed that the lookahead EI is a dynamic program that chooses the EI
maximizer in expectation, considering all possible strategies of the same
budget. Our oracle changes with time and hence the
improvement-based acquisition functions (e.g., EI and lookahead EI) do not suit our problem, since the
appropriate target for a future $t$ is unknown. A practical challenge with
lookahead strategies is the computational tractability, as demonstrated in, for
example, \cite{ginsbourger2010towards} and \cite{osborne2010bayesian}, which is
partially mitigated by various approximation methods,
e.g.,~\cite{wu2019practical, lam2016bayesian, gonzalez2016glasses}.
In this work, we are interested in extending the finite-budget BO of
\cite{osborne2010bayesian} to time-dependent oracles (that is, where the best
action depends on time), in addition to providing a practical solution to the
tractability of the approach when looking ahead more than two steps.

Our main contributions are as
follows: (i) We present the first strategy to solve the finite-budget BO for
time-dependent problems. (ii) We replace the improvement-based acquisition function 
with the GP posterior (payoff) at $T$, with which we are able to obtain
unbiased estimators for our acquisition function and its gradient, which can
then be used in efficient gradient-based optimization methods. (iii) We
improve the tractability of the multi-step lookahead
problem by introducing a \emph{recursive two-step} lookahead strategy 
that looks ahead to $T$ from the current
step $t$; we call our acquisition function Recursive Two-Step Lookahead Expected Payoff
(\texttt{r2LEY}).
(iv) We introduce software that
implements \texttt{r2LEY} for solving time-dependent BO
problems.
(v) Finally, we establish the utility of our method by comparing it against a suite
of existing BO acquisition functions on synthetic and real-world datasets.

\textbf{Related work:} Our approach is similar in spirit to the entropy-based
approaches \cite{hennig2012entropy, hernandez2014predictive, wang2017max} in
the sense that we make decisions that would maximize the expected payoff at $T$
similarly to the way that those approaches make decisions that reduce uncertainty about the global
maximum or maximizer. However, unlike the entropy-based approaches, our approach is able to generate realizations of
analytically closed-form expressions that yield unbiased estimates of our
acquisition function and its gradient. In the context of
time-dependent BO, the only works we are aware of are
\cite{nyikosabayesian} and \cite{bogunovic2016time}. 
We fundamentally differ
from their approach in the following ways: (i) their goal is to track a
time-varying optimum (ours is to predict it at $T$) and (ii) they don't address
the finite-budget constraint as we do. Furthermore, setting $k_t(t_i,
t_j)=(1-\epsilon)^{|t_i - t_j|/2}$ (where $\epsilon$ is the \emph{forgetting}
factor as defined in \cite{bogunovic2016time}), our method is equivalent to
\cite{bogunovic2016time} if we myopically make decisions at each $t$ via UCB.  



\begin{algorithm2e}[H]
\textbf{Given:} \\
\quad $\mcl{D}_n = \lbrace \x_i, y_i \rbrace ~\forall i=1,\hdots,n$, \\
\quad $q$ total budget, \\
\quad schedule $\lbrace t_{n+1},\ldots, t_{q}:= T \rbrace$,\\
\quad and GP hyperparameters $\Omega$ \\
\KwResult{$(\x^*_T, y_T)$}
  \For{$\ell=n+1, \ldots, q$, }{
  Find $\x^*_{\ell} = \underset{\x \in \mcl{X}}{arg\T{max}}~ \alpha(\x, t_{\ell})$ (acquisition function)\\
    Observe $y_{\ell}$ = $f(\x^*_{\ell}, t_{\ell}) + \epsilon_{\ell}$\\ 
    Append $\mcl{D}_{\ell} = \mcl{D}_{\ell-1} \cup \lbrace (\x^*_{\ell}, y_{\ell}) \rbrace$\\ 
    Update GP hyperparameters $\Omega$ \\
 }
 \caption{Generic Bayesian Optimization}
 \label{a:BO}
\end{algorithm2e}
\section{\texttt{r2LEY}: Recursive two-step lookahead expected payoff acquisition function}
Given $n$ observations $\mcl{D}_n$,
let us define a function $I(\x,t)$ which implicitly depends on $Y_n(\x,t)$ and hence is a random
variable. The one-step acquisition function to make the $n+1$th
observation is then defined as
\begin{equation}
    \alpha(\x,t_{n+1}) = \int_{Y_n} I(\x, t_{n+1}) p( Y_n(\x, t_{n+1})|\mcl{D}_n)~dY_n.
    \label{e:acq}
\end{equation}
Note that setting $I(\x, t_{n+1}) = \left(Y_n(\x,t_{n+1}) - \xi(t) \right)^+$
leads to a time-dependent EI; see supplementary material for derivation. In the
general EI setting, the target $\xi$ is the best observed value, which does
not suit time-dependent problems: We seek the maximizer of the
payoff over $\x \in \mcl{X}$ at $t=T$ as opposed to seeking the maximizer over the joint $(\x,t)$ space.
Choosing the target $\xi$ as the maximum of the GP posterior mean---as
in, for instance, \cite{wang2014theoretical}---involves an additional
global optimization step just to determine $\xi$, which could be inaccurate. 

Instead, 
we set
the acquisition function to be the GP posterior at $T$, that is, $I(\x,t_{n+1}) =
Y_n(\x, T)$. Essentially, at each stage, we make a decision that seeks
the best payoff at $T$, our target time horizon. Note that in ~\eqref{e:acq}
and what follows, $Y_n$ refers to the GP posterior conditioned on the past $n$
observations and $\mbb{E}_n$ refers to the expectation taken with respect to
$Y_n$.

\begin{figure}[htb]
    \centering
    \includegraphics[width=\textwidth, trim=0.5cm 4cm 0.5cm 6cm,clip]{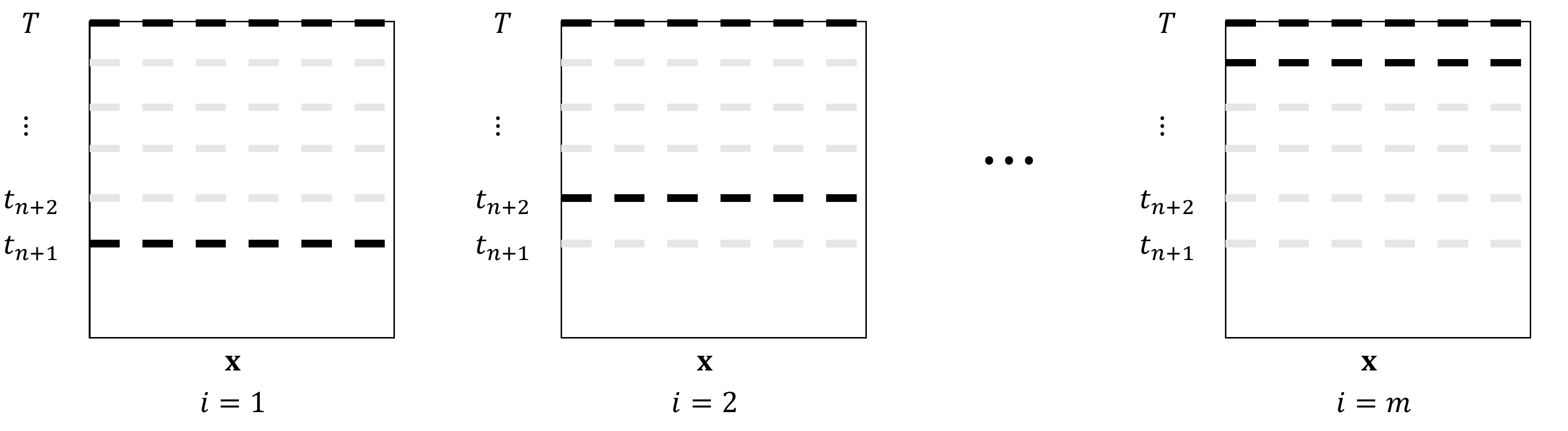}
    \caption{Recursive two-step lookahead method for time-dependent problems.
    Here, we start with $n$ observations up to $t_n$ and have $m$ remaining
    observations to make, where $m+n=q$ (total budget).}
    \label{f:method}
\end{figure}

Following \eqref{e:acq}, at the end of $i=q-1$ rounds of the optimization,
the (final) $q$th observation at $t_{q} := T$ is selected as the
point that maximizes the mean of the GP posterior at time $T$. That is,
\begin{equation}
    \begin{split}
        \x^*_{q} =& \underset{\x_{q} \in \mcl{X}}{arg\T{max}} \int_{Y_{q-1}}Y_{q-1}(\x_{q},T)p(Y_{q-1}|\mcl{D}_{q-1})~dY_{q-1} \\
         =& \underset{\x_{q} \in \mcl{X}}{arg\T{max}}~\mbb{E}_{q-1} Y_{q-1}(\x_{q},T) \\
         =& \underset{\x_{q} \in \mcl{X}}{arg\T{max}} ~\mu_{q-1}(\x_{q}, T),
    \end{split}
\end{equation}
where $p(Y_{q-1}|\mcl{D}_{q-1})$ is the density of $Y_{q-1}$. Going back
by one time-step to after $q-2$ observations of the oracle have been performed, the
point $\x^*_{q-1}$ is selected as follows.\footnote{While $\x_{q}^*$
is an estimate of the maximizer of $f$ at $t_{q}$, $\x^*_{q-1}$ is not an
estimate of the maximizer of $f$ at $t_{q-1}$. We use this notation to
disambiguate the solution to the optimization problem in \eqref{eq:nm1} from
its decision variable $\x_{q-1}$.} 
Note
that here, the $q-1$th observation $y_{q-1}$ is unknown and hence is drawn
from the distribution $\mcl{N}(\mu_{q-2}(\x, t_{q-1}), \sigma^2_{q-2}(\x,
t_{q-1}))$ given whatever $\x$ we choose to observe at $t_{q-1}$.
\begin{equation}\label{eq:nm1}
    \begin{split}
    \x^*_{q-1} = \underset{\x_{q-1}\in\mcl{X}}{arg\T{max}}~\int_{y_{q-1}} &\left[ \underset{\x_{q} \in \mcl{X}}{\T{max}} \mbb{E}_{q-1} Y_{q-1}(\x_{q},T)|\mcl{D}_{q-2}, y_{q-1} \right] \\
    & p(y_{q-1}|\x_{q-1}, \mcl{D}_{q-2})~dy_{q-1},
    \end{split}
\end{equation}
where, in \eqref{eq:nm1}, the inner maximization is with respect to
$\x_{q}$ and the outer maximization is with respect to $\x_{q-1}$. The above
equation is written concisely as
\begin{equation}
    \x^*_{q-1} = \underset{\x_{q-1} \in \mcl{X}}{arg\T{max}}~\mbb{E}_{q-2} \left[ \underset{\x_{q} \in \mcl{X}}{\T{max}} \mbb{E}_{q-1} Y_{q-1}|\mcl{D}_{q-2}, y_{q-1} \right].
    \label{e:two-step}
\end{equation}
Similarly, we could extend \eqref{e:two-step} to look ahead three steps to
determine $\x^*_{q-2}$, which, however, would involve the marginalization of
$\x_{q-1}$:
\begin{equation}
    \begin{split}
    \x^*_{q-2} = \underset{\x_{q-2}\in\mcl{X}}{arg\T{max}}~\int_{y_{q-2}} \int_{\x_{q-1}} &\left[ \underset{\x_{q} \in \mcl{X}}{\T{max}} \mbb{E}_{q-1} Y_{q-1}(\x_{q},T)|\mcl{D}_{q-2}, y_{q-1} \right] \\
    & p(y_{q-1}|\x_{q-1}, \mcl{D}_{q-2})\\
    & p(\x_{q-1}|\mcl{D}_{q-2})~d\x_{q-1}~dy_{q-1},
    \end{split}
\end{equation}
where $p(\x_{q-1}|\mcl{D}_{q-2}) = p(\x_{q-1})$ can be set as
$\mcl{U}(lb, ub)$, where $lb$ and $ub$ are the lower and upper
bounds of the domain $\mcl{X}$.
This additional step of marginalizing
$\x_{q-1}$ when moving from the two- to three-step lookahead makes the method
computationally intractable, particularly as $d$ increases. Using the same
principle, one can define a generalized $m$-step lookahead to select $\x^*_{q-m}$ by appropriately
marginalizing all the remaining $(\x_i, y_i)$, $\forall i=q-m+1,\ldots,q$. 
Instead, in this work, we recursively apply the two-step lookahead approach
~\eqref{e:two-step} to multi-step time-dependent BO problems (see
Figure~\ref{f:method}) and show how to efficiently optimize the two-step
acquisition function defined as follows.

First, we write our two-step lookahead payoff acquisition function to select
the $j+1$th point---with $m=q-j$ observations remaining---as follows:
\begin{equation}
    \begin{split}
    \alpha_{2LEY}(\x_{j+1}) =& \mbb{E}_j \left[ \underset{\x_{q} \in \mcl{X}}{\T{max}}~ \mu_{j+1}(\x_{q}, T) | y_{j+1}, \x_{j+1} \right] \\
    =& \mbb{E}_j \left[ \underset{\x_{q} \in \mcl{X}}{\T{max}}~ g(\x_{q}, \x_{j+1}, T) \right],
    \end{split}
    \label{e:2LEY}
\end{equation}
where $g(\x_{q}, \x_{j+1}, T) = \mu_{j+1}(\x_{q}, T) | y_{j+1}, \x_{j+1}$ and
$y_{j+1}$ is a draw from $\mcl{N} \left(\mu_j(\x, t_{j+1}), \sigma_j^2(\x,
t_{j+1}) \right)$. Let $g^*(\x_{j+1}, T) = g(\x^*_{q}, T, \x_{j+1})$, where
$\x^*_{q}$ is a maximizer of $\mu_{j+1}(\x_{q}, T)$. Then the above equation
can be written
\begin{equation}
    \alpha_{2LEY}(\x_{j+1}) = \mbb{E}_j \left[ g^*(\x_{j+1}, T) \right].
    \label{e:2LEY_2}
\end{equation}
Note that in \eqref{e:2LEY_2}, it is implicit that $\alpha_{2LEY}(\x_{j+1})$ is
the acquisition function to choose a point at $t_{j+1}$. The dependence of the
right-hand side of \eqref{e:2LEY_2} and the second line of \eqref{e:2LEY} on
$\x_{j+1}$ can
be seen by realizing that
\begin{equation}
    g^*(\x_{j+1}, T) = \mu_{j+1}(\x^*_{q}, T) = \mbf{k}_{j+1}^\top \mbf{K}^{-1}_{j+1} \mbf{y}_{j+1},
\end{equation}
where $ \mbf{k}_{j+1} = [ k\left((\x^*_{q},T), (\x_1, t_1)\right),
\ldots,
k\left((\x^*_{q},T), (\x_j, t_j)\right), k\left((\x^*_{q},T), (\x_{j+1},t_{j+1})\right)]^\top$. The expectation in \eqref{e:2LEY_2} is not available in closed form, but can be approximated as
\begin{equation}
    \alpha_{2LEY}(\x_{j+1}) = \mbb{E}_j \left[ g^*(\x_{j+1}, T) \right] \approx \frac{1}{M} \sum_{i=1}^M  g^*(\x_{j+1}, T) | y^i_{j+1}
    \label{e:2LEY_MC}
\end{equation}
The gradient of $g^*(\x_{j+1}, T)$ with respect to $\x_{j+1}$ is given by
\begin{equation}
    \begin{split}
        \nabla g^*(\x_{j+1}, T) = &\f{\partial \mbf{k}_{j+1}^\top}{\partial \x_{j+1}} \mbf{K}^{-1}_{j+1} + \mbf{k}_{j+1}^\top \f{\partial \mbf{K}^{-1}_{j+1}}{\partial \x_{j+1}} \\
        = &\f{\partial \mbf{k}_{j+1}^\top}{\partial \x_{j+1}} \mbf{K}^{-1}_{j+1} +
        \mbf{k}_{j+1}^\top~\mbf{K}^{-1}_{j+1} \f{\partial \mbf{K}_{j+1}}{\partial \x_{j+1}}\mbf{K}^{-1}_{j+1},
    \end{split}
    \label{e:2LEY-grad}
\end{equation}
where $ \f{\partial \mbf{K}_{j+1}}{\partial \x_{j+1}}$ is a matrix of
elementwise derivatives and the second line follows from a well-known lemma on the derivative of matrix inverse~\cite[p.~201--202]{rasmussen:williams:2006}.
With a continuously differentiable kernel $k(\cdot, \cdot)$ we state that
\begin{equation}
    \nabla \alpha_{2LEY}(\x_{j+1}) = \nabla \mbb{E}_j \left[ g^*(\x_{j+1}, T) \right] = \mbb{E}_j \left[\nabla g^*(\x_{j+1}, T) \right],
    \label{e:2LEY_gr_ex_interchange}
\end{equation}
where the interchange of the gradient and expectation operators is via Theorem~\ref{thm:grad_exp}. Then the gradient is approximated as 
\begin{equation}
    \nabla \alpha_{2LEY}(\x_{j+1}) = \mbb{E}_j \left[\nabla g^*(\x_{j+1}, T) \right] \approx \frac{1}{M} \sum_{i=1}^M g^*(\x_{j+1}, T) | y^i_{j+1}.
    \label{e:2LEY_grad_MC}
\end{equation}
Furthermore, it can
be shown (see Theorem~\ref{thm:grad_exp}) that \eqref{e:2LEY_MC} and \eqref{e:2LEY_grad_MC} are unbiased estimators for $\alpha_{2LEY}$ and $\nabla
\alpha_{2LEY}$, respectively. This gradient can then be used in a gradient-based optimizer to
maximize our two-step lookahead acquisition function in an efficient manner.
The overall algorithm is given by inserting the acquisition function obtained
via Algorithm~\ref{a:2LEY1} into line 7 of Algorithm~\ref{a:BO}.

\begin{algorithm2e}[H]
\SetAlgoLined
\textbf{Given:} 
point $\x_{j+1}$,
data $\mcl{D}_{j} = \{ (\x_i, t_i), y_i \} \forall i=1,\hdots,j$,
and (simulation) budget $M$\\
\KwResult{estimate of $\alpha_{2LEY}$ and $\nabla \alpha_{2LEY}$}
\For{$i = 1,\ldots,M$}{
    $y_{j+1}^i \sim \mcl{N}(\mu_{j}, \sigma^2_{j})$ \hspace*{\fill} (simulated $y_{j+1}$)\\
    $\mcl{D}_{j+1}^i = \mcl{D}_{j} \bigcup \{(\x_{j+1}, y^i_{j+1}) \}$ \hspace*{\fill} (update data with simulation) \\
    $Y^i_{j+1}|\mcl{D}_{j+1}^i \sim \mcl{GP}(\mu^i_{j+1}, \sigma^{2~i}_{j+1})$ \hspace*{\fill} (update GP posterior)\\
    $\x_q^{*i} = \underset{\x \in \mcl{X}}{arg\T{max}}~\mu^i_{j+1}(\x, \x_{j+1}, T)$ \hspace*{\fill} (maximizer of payoff at $T$)\\
    $\nu^i = \mu^{i*}_{j+1}(\x_{j+1}, T)$ \hspace*{\fill} (simulated acquisition function) \\
    $d\nu^i = \f{\partial \mu^{i*}_{j+1}(\x_{j+1}, T)}{\partial \x_{j+1}}$ \hspace*{\fill} (simulated acquisition function gradient) \\
}  
$ \alpha_{2LEY} = \mbb{E}_j \left[ g^*(\x_{j+1}, T) \right] \approx \frac{1}{M} \sum_{i=1}^M \nu^i$ \\
$ \nabla \alpha_{2LEY}(\x_{j+1}) = \mbb{E}_j \left[ \partial g^*(\x_{j+1}, T)/ \partial \x_{j+1} \right] \approx \frac{1}{M} \sum_{i=1}^M d\nu^i$
 \caption{Monte Carlo Approximation of \eqref{e:2LEY} and \eqref{e:2LEY_gr_ex_interchange}} 
 \label{a:2LEY1}
\end{algorithm2e}
\section{Theoretical Properties}
\begin{theorem}[Interchange of gradient and expectation operators] Let us write
  $g^*(\x_{j+1}, T)|y_{j+1}$ as $g^*(\x_{j+1}, T, y_{j+1})$ and, let $\mcl{Y}$
  be the support of $y_{j+1}$, whose density is $p(y_{j+1})$, and suppose the
  domain of $\x_{j+1}$ $\mcl{X}$ is an open set. Let $g^*(\x_{j+1}, T, y_{j+1})
  p(y_{j+1})$ and $\partial g^*(\x_{j+1}, T, y_{j+1})/\partial \x_{j+1}
  p(y_{j+1})$ be continuous on $\mcl{X} \times \mcl{Y}$. Suppose that there
  exist nonnegative functions $q_0(y_{j+1})$ and $q_1(y_{j+1})$ such that
  $|g^*(\x_{j+1}, T, y_{n+1}) p(y_{j+1})| \leq q_0(y_{j+1})$ and $\|
  \frac{\partial g^*}{\partial \x_{j+1}}  p(y_{j+1})\| \leq q_1(y_{j+1})$ for
  all $(\x_{j+1}, y_{j+1}) \in \mcl{X} \times \mcl{Y}$, where
  $\int_{\mcl{Y}}q_1(y_{j+1}) dy_{j+1} < \infty$ and
  $\int_{\mcl{Y}}q_2(y_{j+1}) dy_{j+1} < \infty$. Then,
\[ \nabla \mbb{E}_j \left[ g^*(\x_{j+1}, T, y_{j+1}) \right] = \mbb{E}_j \left[\nabla g^*(\x_{j+1}, T, y_{j+1}) \right]\]

Since $\nabla \alpha_{2LEY}(\x_{j+1}) = \mbb{E}_j \left[\nabla g^*(\x_{j+1}, T,
y_{j+1}) \right]$, a realization of $\nabla g^*(\x_{j+1}, T, y_{j+1})$ yields
an unbiased estimate of the true gradient.
\label{thm:grad_exp}
\end{theorem}
See supplementary material for proof.

\section{Experiments}
The synthetic and real-world test data used in the experiments are described as
follows. We introduce new synthetic test functions that can serve to benchmark
general time-varying stochastic optimization problems.

\begin{figure}[ht]
    \centering
    \begin{subfigure}[htb!]{.32\textwidth}
        \includegraphics[width=\textwidth]{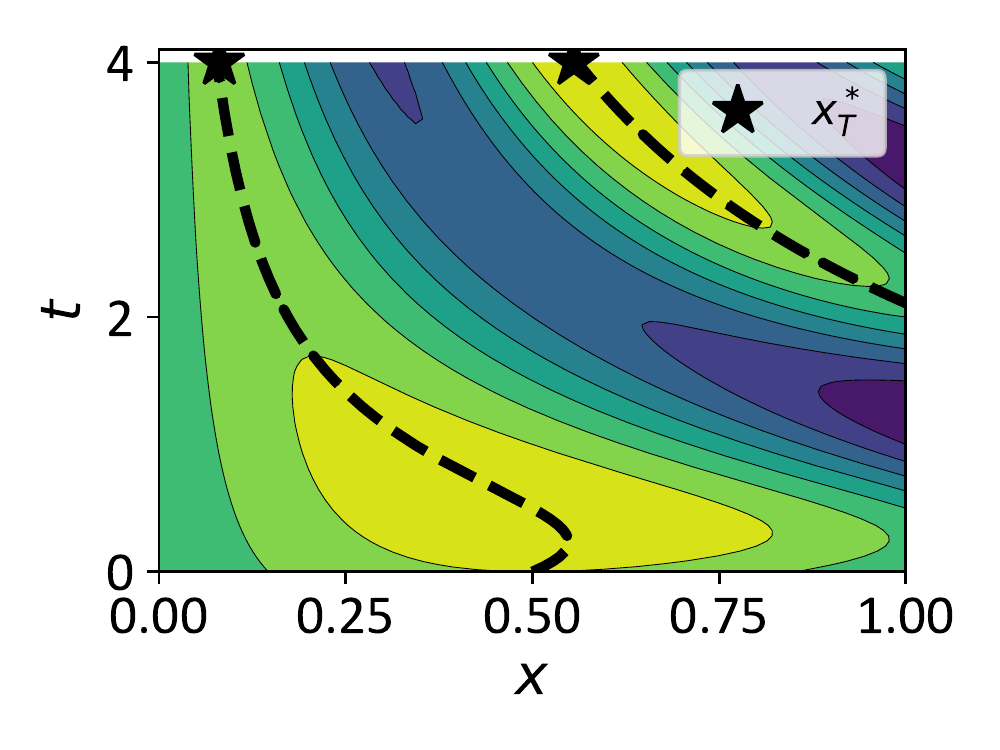}
        \caption{Quadratic-b}
        \label{sf:qb}
    \end{subfigure}
    \begin{subfigure}[htb!]{.32\textwidth}
        \includegraphics[width=\textwidth]{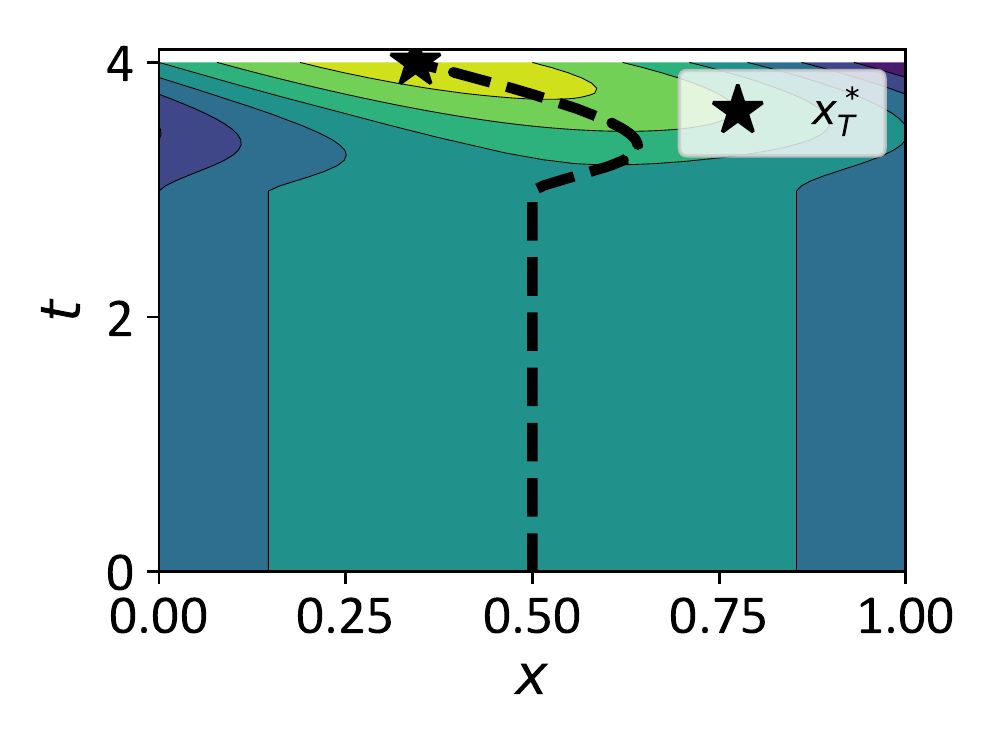}
        \caption{Quadratic-c}
        \label{sf:qc}
    \end{subfigure}
    \begin{subfigure}[htb!]{.32\textwidth}
        \includegraphics[width=\textwidth]{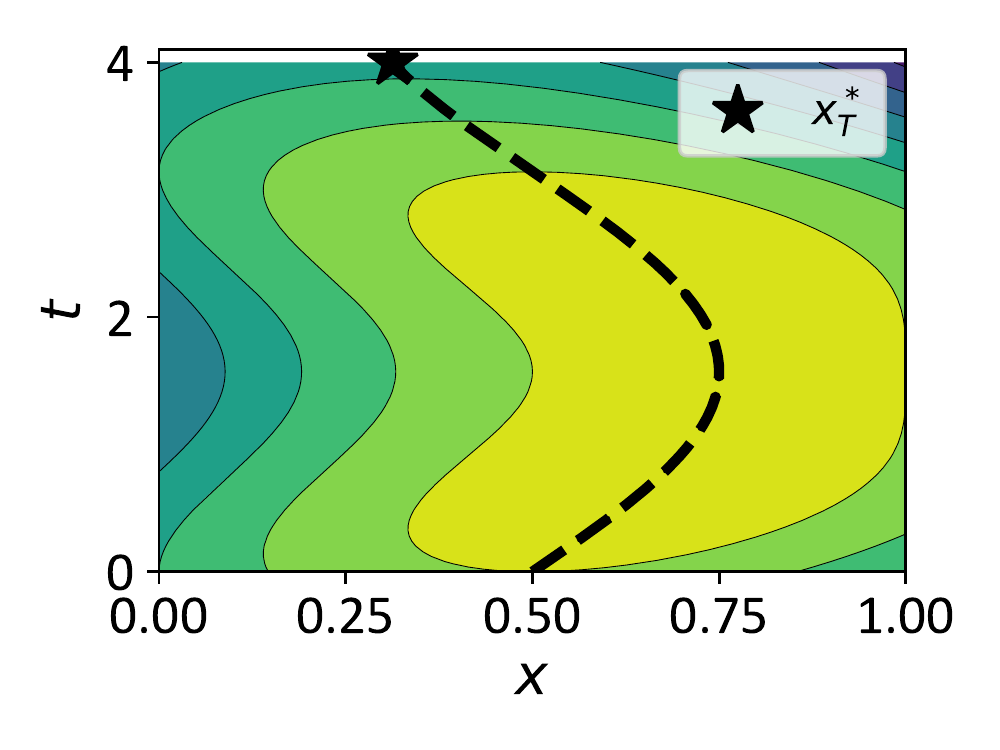}
        \caption{Quadratic-d}
        \label{sf:qd}
    \end{subfigure}

    \caption{Contours that visualize the one-dimensional ($\mcl{X} \subset
    \mbb{R}$) time-dependent payoff function $f(\x,t)$. The dashed lines
    trace the location of $\x^*$ at every $t$. We are
    interested in the maximizer of the payoff at $T$, $\x^*_T$. Each payoff
    function is generated by adding a nonlinear component (dependent on $\x$ at
    $t$) to a simple one-dimensional quadratic function; see supplementary
    material for details.} 
    \label{f:1d_quad}
\end{figure}

\textbf{Synthetic one-dimensional test functions}: Each synthetic test function
takes the form $f(\x, t) = f_{\x}(\x) + f_{\x t}(\x, t)$, where $f_{\x t}$ adds
context dependence. For the one-dimensional cases, $f_{\x}(\x) =
-4\times(\x-0.5)^2$ and $f_{\x t}$ is chosen to generate varying trajectories
for $\x^*$ at each $t$. The specific details are provided in the supplementary
material and the contour plots of the payoff function are shown in
Fig.~\ref{f:1d_quad}. Notice that Quadratic c \& d have exactly one local
maximizer at every $t$ whereas Quadratic b can be multimodal at specific
$t$'s. 

\textbf{Synthetic higher-dimensional test functions}: For dimensions $d>1$,
various test functions are chosen for $f_{\x}(\x)$, whereas $f_{\x t}$ is set
to be that of Quadratic-d. The
synthetic functions (for $f_{\x}(\x)$) chosen are the Griewank
(2d), Hartmann (3d), Hartmann (6d), Levy (8d), and Syblinski-Tang (10d); see
\cite{simulationlib} for details of each of these test functions.

\textbf{Intel temperature sensor data}: The following data are from sensor
measurements from the Intel Berkeley Research Center, where 47 sensors are
distributed in the building to measure various ambient air properties. We are
specifically interested in the temperature data that are available as a time
series for each sensor. 
We use the first
50 hours to train our model and predict the maximum temperature at $T=100$.

\textbf{SARCOS robot arm data}: These data relates to an inverse dynamics
problem for a SARCOS anthropomorphic robot arm with seven degrees of freedom. The
original data have a 21-dimensional input space (7 joint positions, 7 joint
velocities, 7 joint accelerations) and a 7-dimensional output space (joint
torques). We use the first 7 of the 21 input
variables as the $\x$, the $8$th input variable as $t$ and the first of the seven torques as output. The input $t$ is scaled to range from $t=0$ to $T=4$ and the whole data is sorted (in ascending order) with respect to $t$. Data up to $t=2$ are used as initial
samples and the horizon is set as $T=4$.

We start each experiment with $n=(d+1)\times 20$ initial samples for test functions upto $d=6$ and $n=(d+1)\times 10$ for higher dimensions and, run our
algorithm for $m=10$ additional steps at fixed, equal time-steps. All initial
samples span $0\leq t\leq 2$ with one sample per time at equally spaced
time-steps and the horizon is set at $T=4$, unless otherwise mentioned.
$M=5000$ is set for all experiments and each experiment is repeated 20 times.
Our metric for comparison is the $\mathtt{log_{10}}$ normalized-simple-regret
at $T$ defined as $\log_{10}~ \frac{f_{max} - f(\x^*_T)}{f_{max} -
f_{min}}$, where $f_{max}$ and $f_{min}$ are the maximum and minimum values for
$f(\x, T)~\forall \x \in \mcl{X}$, and $\x^*_T$ is the point selected by the
algorithm at $T$. Both the mean (with standard errors) and median (with quantiles) values of the metric are compared
across all experiments.

We compare our method against the most widely used myopic approaches in BO,
namely the EI, PI, and UCB, which are modified appropriately for time-dependent
problems (see supplementary material).
\texttt{EImumax} and \texttt{PImumax} are the EI and PI,
respectively with the target set as the maximum of the GP posterior mean at the
current step, that is, $\mu_n(\x_{n+1}, t_{n+1})$, and the confidence parameter for
\texttt{UCB} is set to $2$. Additionally, we compare against a strategy that
selects points uniformly at random from $\mcl \sim \mcl{U}(lb, ub)$, where $lb$
and $ub$ are the lower and upper bounds of a hypercube. Additionally, we
include \texttt{R-EI}, which is essentially \texttt{Random} except for the last
evaluation, which is selected per \texttt{EImumax}.

\subsection{Discussion}
The metrics for all the numerical experiments, computed from 20 independent
replications of each algorithm, are tabulated in Table~\ref{tab:results}.
It is observed that the proposed \texttt{r2LEY} approach outperforms
the other methods in terms of either the best average regret or the
best worst-case regret. See the supplementary material for a more visual
representation of these results.

\begin{figure}[ht]
    \centering
    \begin{subfigure}[htb!]{\textwidth}
        \includegraphics[width=0.2\textwidth]{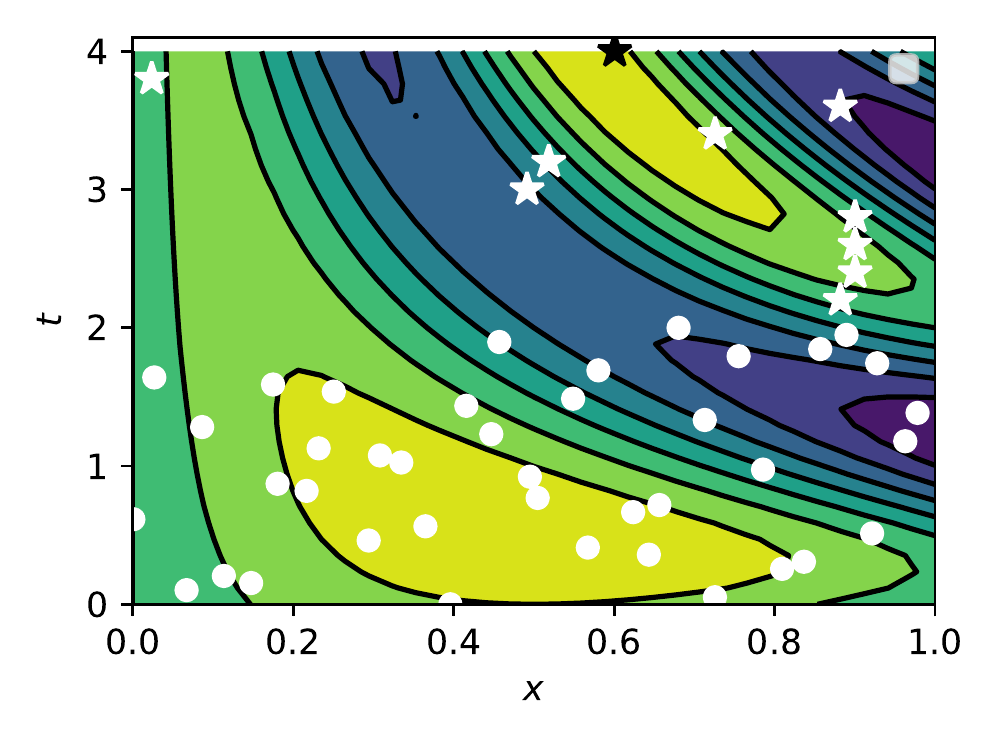}%
        \hfill
        \includegraphics[width=0.2\textwidth]{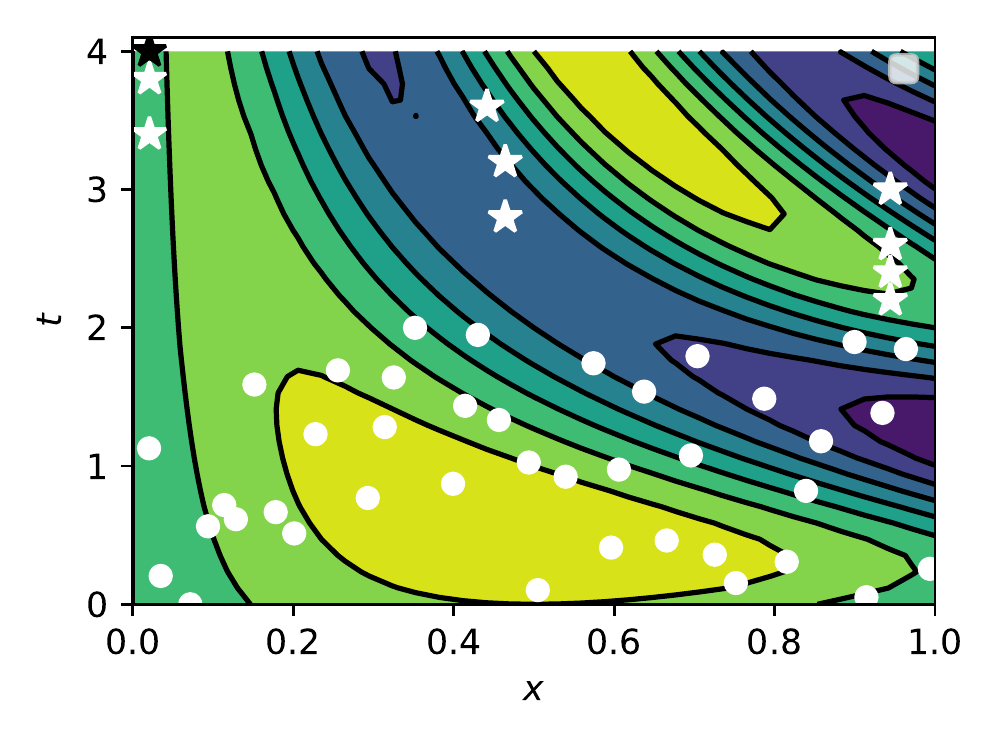}%
        \hfill
        \includegraphics[width=0.2\textwidth]{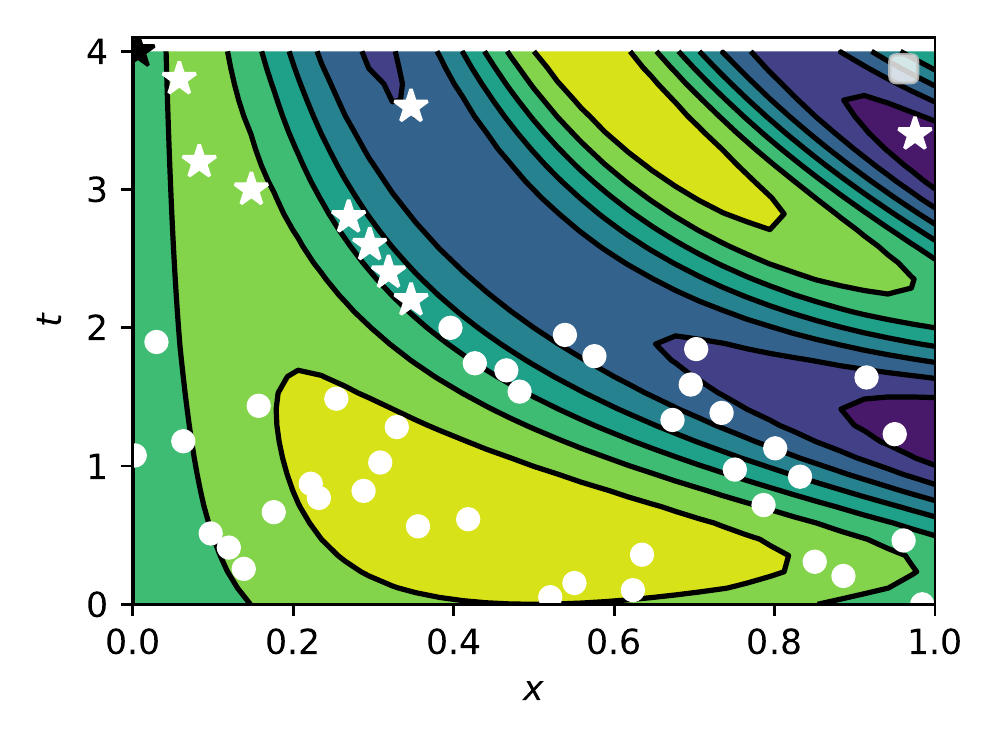}%
        \hfill
        \includegraphics[width=0.2\textwidth]{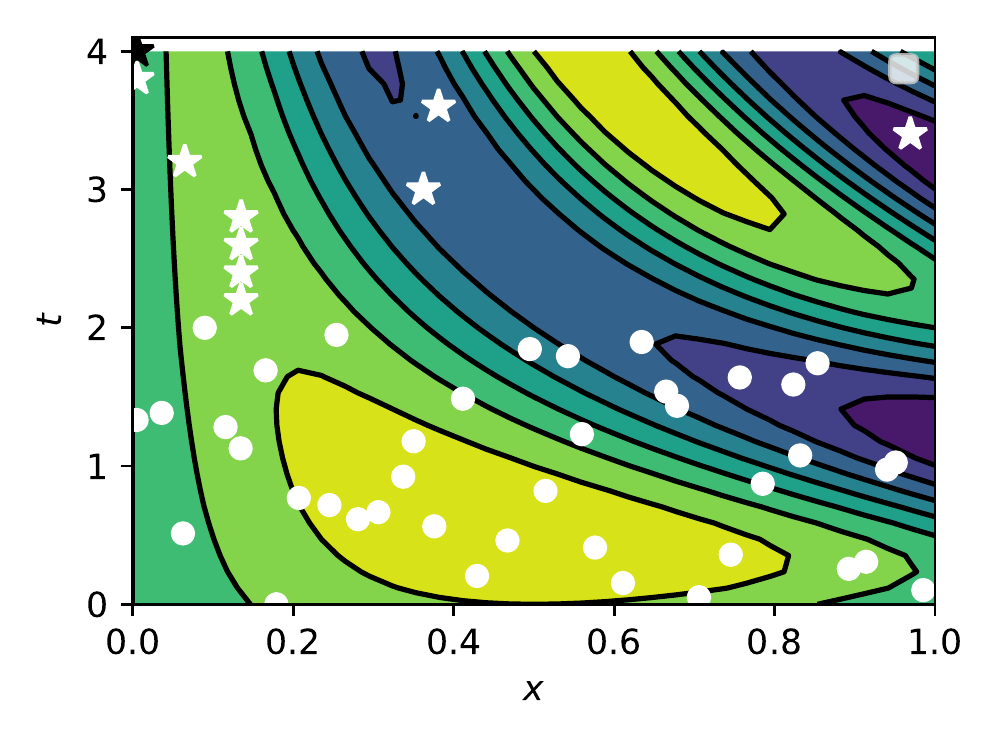}%
        \hfill
        \includegraphics[width=0.2\textwidth]{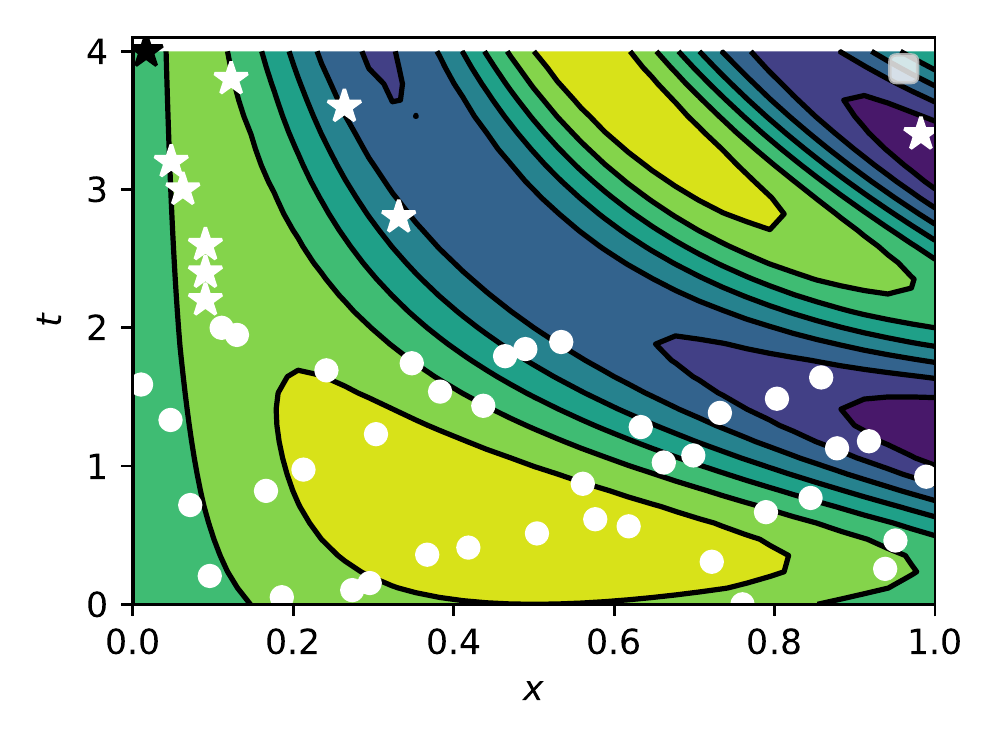}%
        \caption{Quadratic-b}
        \label{sf:qb}
    \end{subfigure}\\
    \begin{subfigure}[htb!]{\textwidth}
        \includegraphics[width=0.2\textwidth]{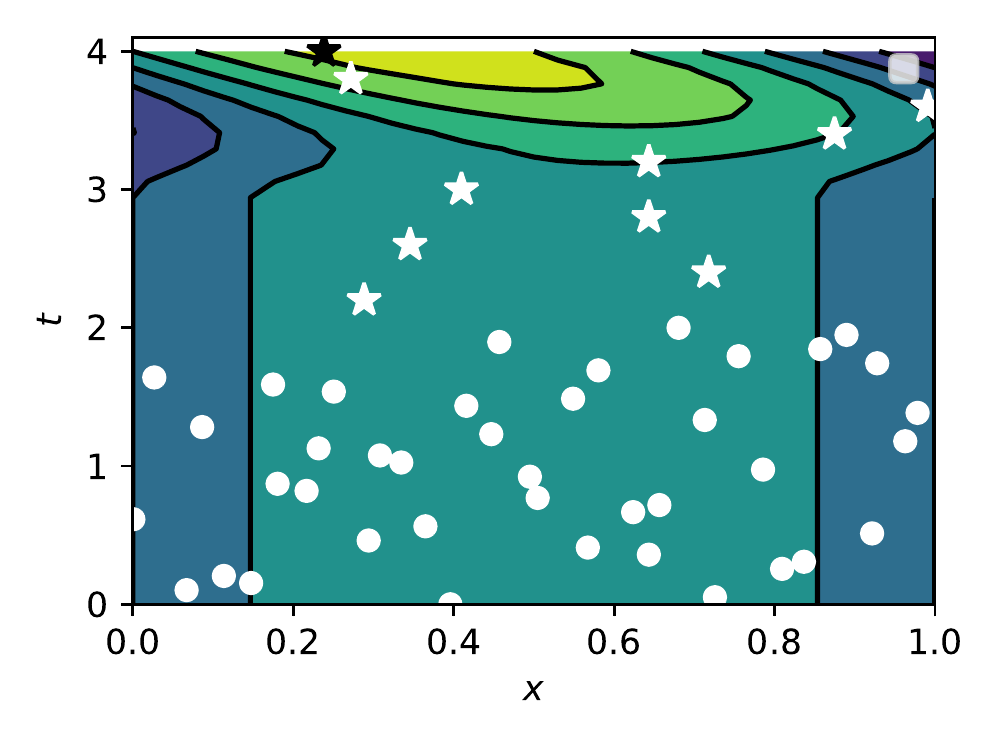}%
        \hfill
        \includegraphics[width=0.2\textwidth]{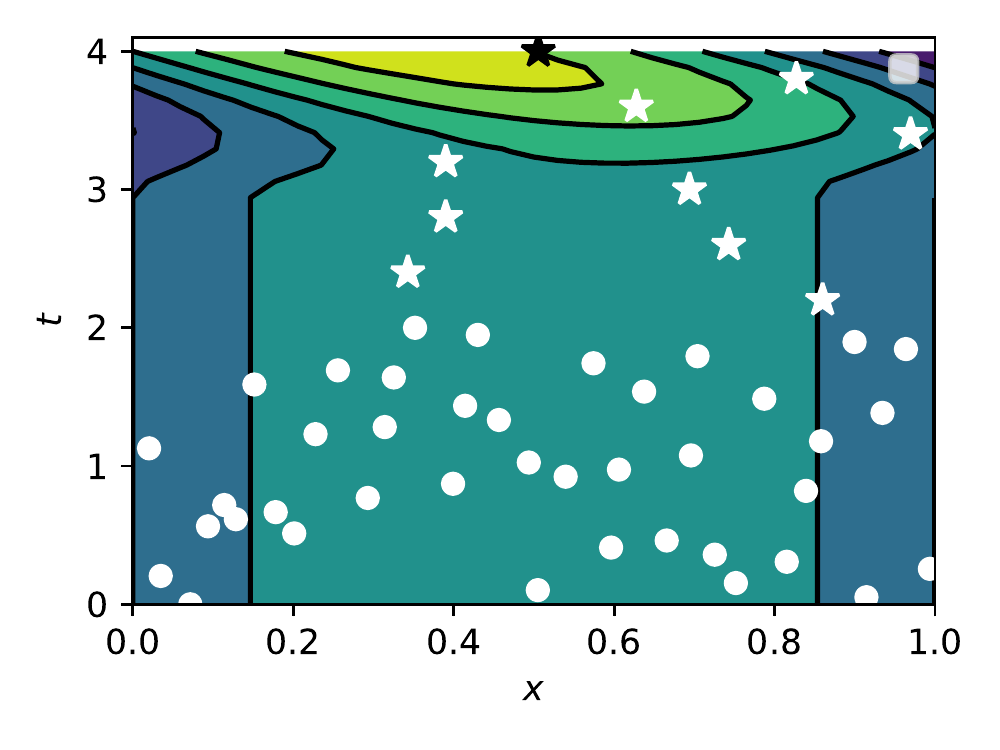}%
        \hfill
        \includegraphics[width=0.2\textwidth]{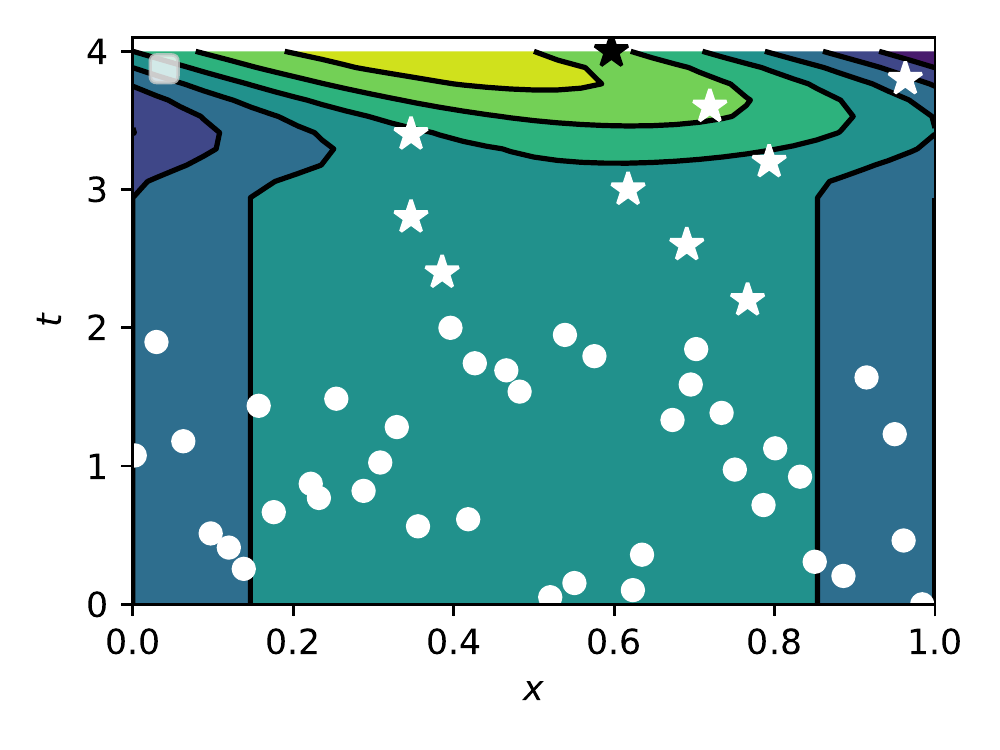}%
        \hfill
        \includegraphics[width=0.2\textwidth]{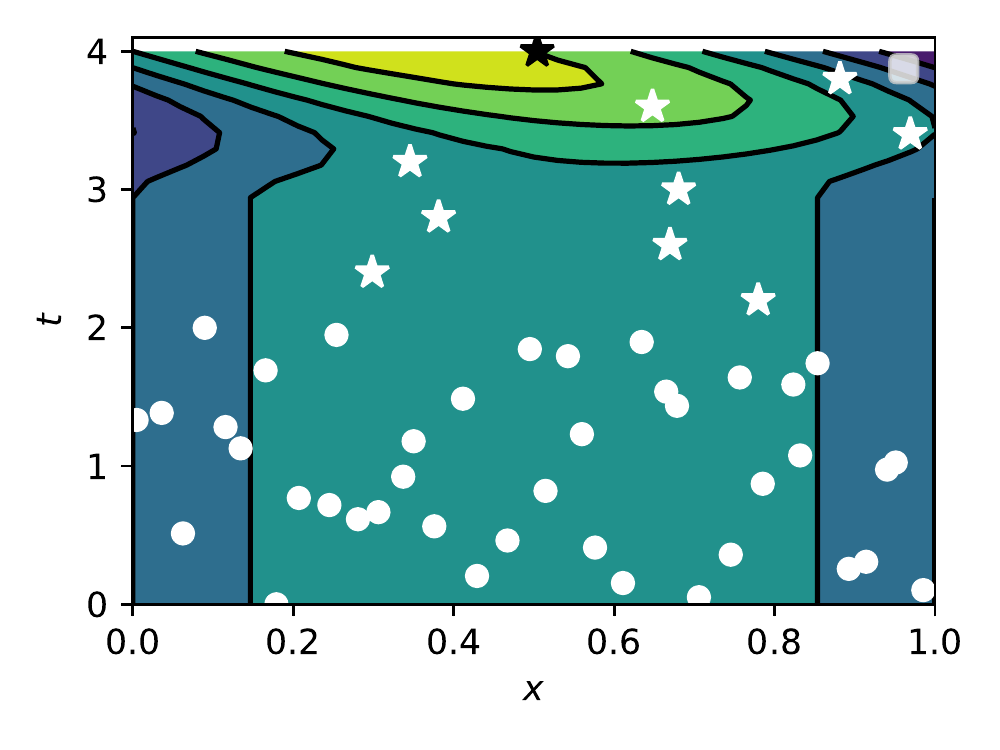}%
        \hfill
        \includegraphics[width=0.2\textwidth]{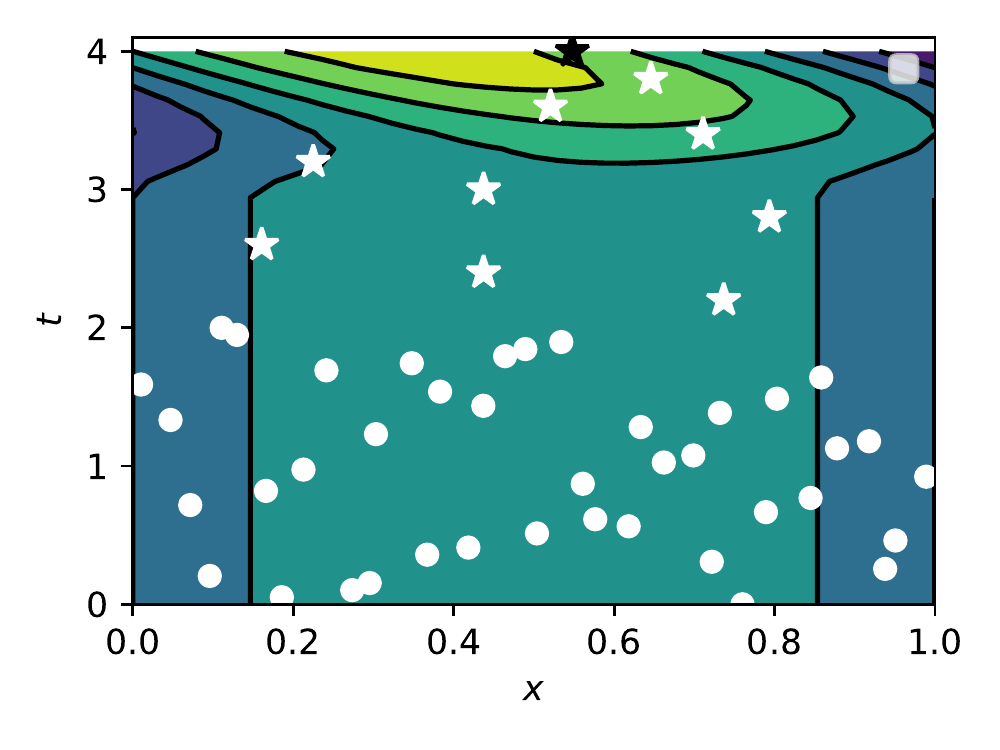}%
        \caption{Quadratic-c}
        \label{sf:qc}
    \end{subfigure}\\
    \begin{subfigure}[htb!]{\textwidth}
        \includegraphics[width=0.2\textwidth]{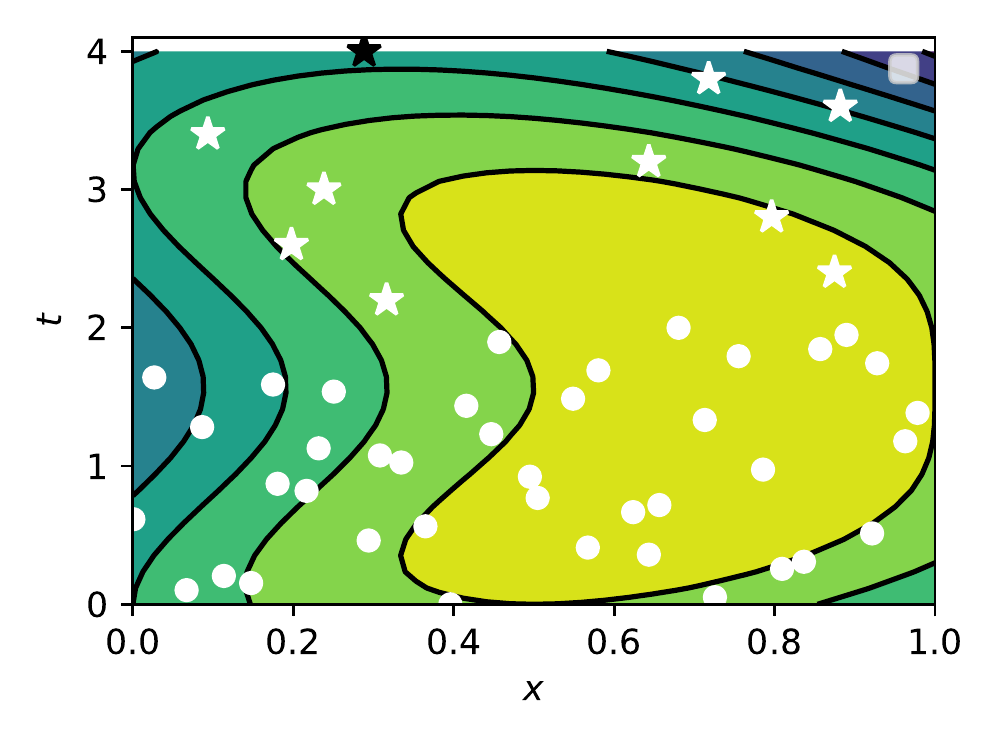}%
        \hfill
        \includegraphics[width=0.2\textwidth]{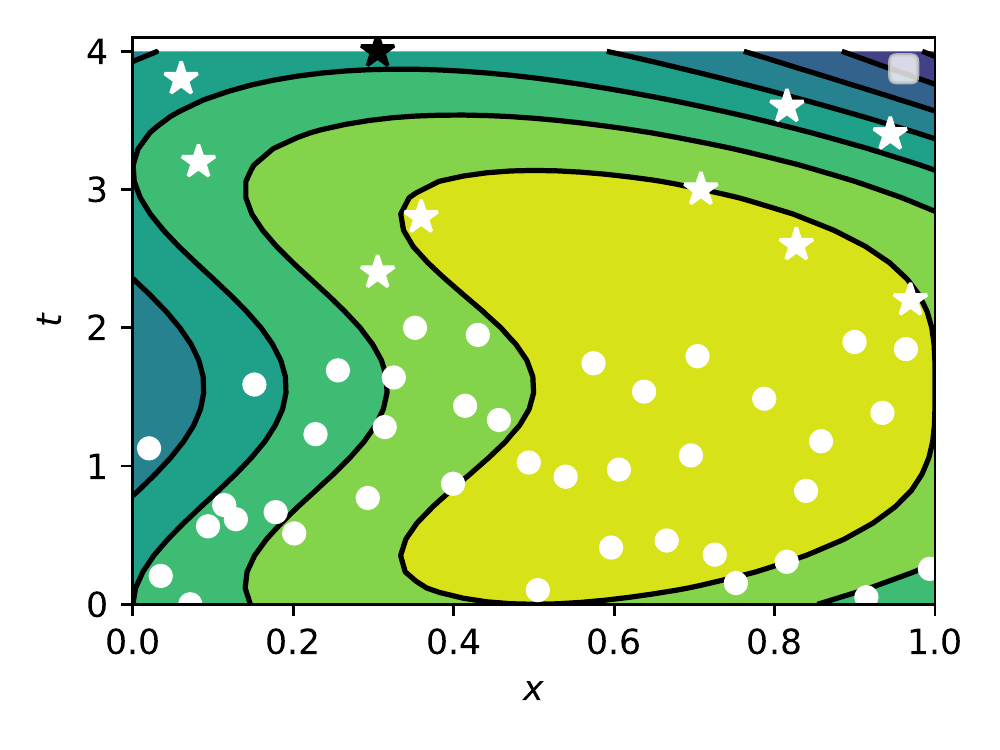}%
        \hfill
        \includegraphics[width=0.2\textwidth]{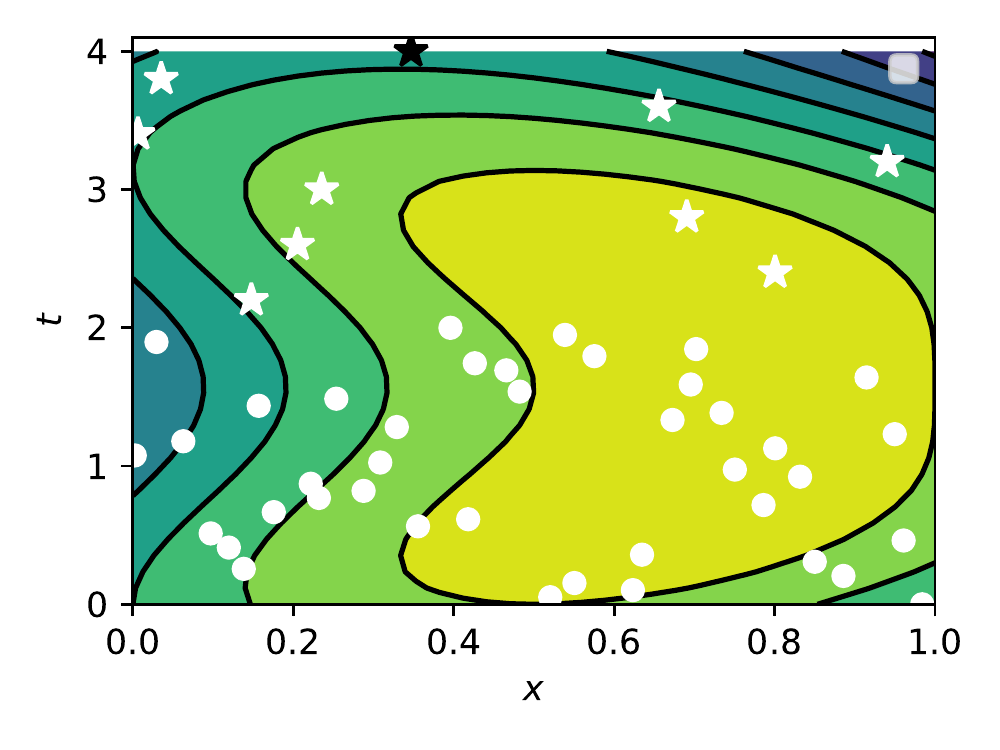}%
        \hfill
        \includegraphics[width=0.2\textwidth]{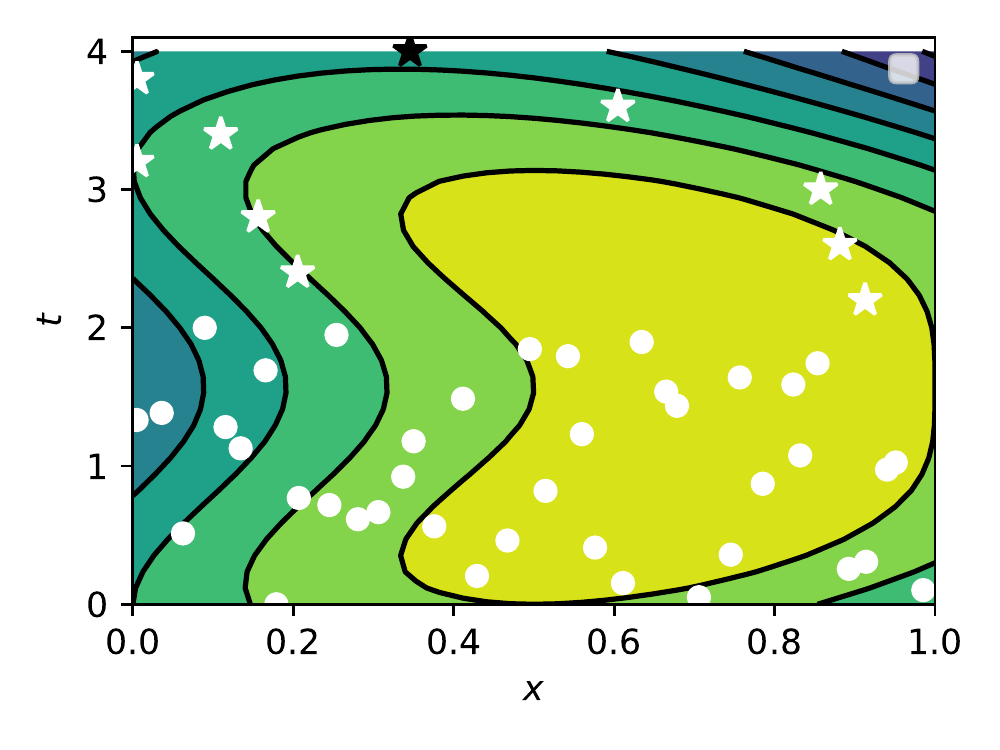}%
        \hfill
        \includegraphics[width=0.2\textwidth]{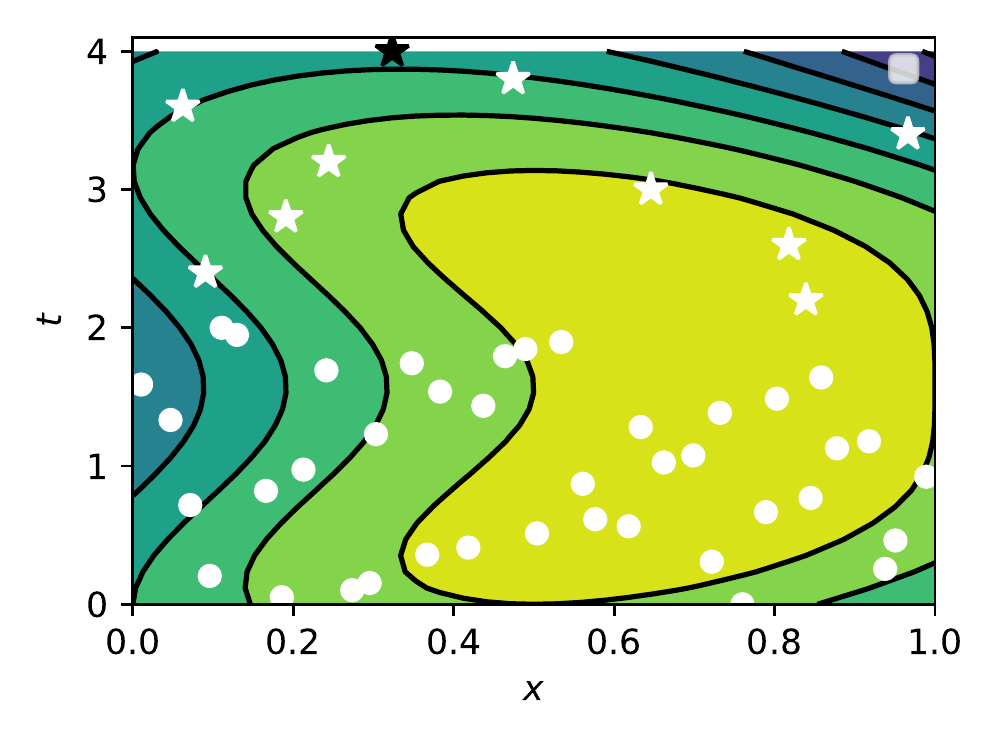}%
        \caption{Quadratic-d}
        \label{sf:qd}
    \end{subfigure}\\
 
    \caption{Performance of the proposed approach on the 1d problems for 5
    repetitions (left to right). Circles are initial samples, stars are points
    selected by \texttt{r2LEY} and the black star is the final point $\x^*_T$.
    Contours represent the \emph{true} noise-free payoff function for
    reference.} 
    \label{f:1d_quad_result}
\end{figure}

The optimizer histories for the $d=1$ (quadratic) test cases are shown in
Figure~\ref{f:1d_quad_result}, where each row represents a unique test function
and each column represents one replication of \texttt{r2LEY}. In these figures, circles
are the starting points and stars are the points chosen via \texttt{r2LEY},
and the black star is the final point at $T$. The contours in the background
correspond to the true $f(\x, t)$. The points placed away from the yellow
regions of the contours are indicative of the property of the method that it
chooses locations of low immediate payoff to ultimately maximize the expected payoff at
$T$. This property indirectly leads the algorithm to favor exploration away from $T$
and exploitation at $T$, suiting our goal of maximizing the payoff solely at
$T$. Among the three test functions in Figure~\ref{f:1d_quad_result}, the
Quadratic-b is particularly challenging because of the sharp
gradients in the $t$ direction; notice the sudden change in the location of the
global maximum (around $t=2$).
Despite these challenges, the proposed approach is able to select the best
payoff at $T$ better than competing methods. In the case of Quadratic-b, the
final point (black star) was always in the neighborhood of one of the two local
maxima at $T$. 
\subsection{Computational Details}
The benefits of the \texttt{r2LEY} are somewhat offset by the high computational
costs relative to the other methods compared in this work. This cost mainly
comes from the \texttt{for} loop in Algorithm~\ref{a:2LEY1}, which includes
$M=5000$ computations of (1) GP posterior update with an additional
observation, (2) GP posterior mean maximization and (3) the gradient
computation of the acquisition function. However, these costs are mitigated in this
work in the following ways. First, the covariance matrix inverse
$\mbf{K}^{-1}_n$ is computed only once and stored, after which, with every new
observation, updates to $\mbf{K}^{-1}_n$ are made to get $\mbf{K}^{-1}_{n+1}$
via the lemma in \cite[p.~77]{press2002numerical}. This procedure reduces the
computation from $\mcl{O}(n+1)^3$ to $\mcl{O}(n+1)^2$. Second, owing to the
independence of each of the $M$ computations in Algorithm~\ref{a:2LEY1}, they
are executed in parallel; in this regard we present a parallel implementation
of \texttt{r2LEY} to be available upon publication. Third, the gradients of
the acquisition function in \texttt{r2LEY} can be directly specified when
analytically known. In other cases, the implementation of \texttt{r2LEY} takes
advantage of automatic differentiation to compute gradients. Fourth, the
maximization of the GP posterior mean is similarly solved efficiently with
multistart gradient-based optimization. Finally, the GP hyperparameter
updates are removed from the \texttt{for} loop and instead done only once per
each of the $m$ time-steps from $t_n$ to $T$. With 72 processors, $M=5000$ and
$m=1$, the current implementation  takes approximately 3 minutes of wall-clock
time to execute one instance of \texttt{r2LEY}. 

\section{Conclusions}
Optimization of a time-dependent expensive  stochastic oracle finds application
in a variety of fields including quantum computing, finance, power \& energy
systems and aerospace engineering. In this work, we specifically address
problems where the best decision at a finite time horizon is of interest, given
a finite budget of evaluations. Our proposed approach extends the lookahead
approach in BO to time-dependent problems. We overcome the tractability issue
of multi-step lookahead approaches by recursively solving a two-step lookahead
problem where, at each step, we look ahead at the desired time horizon.
Furthermore, we solve the optimization of our acquisition function efficiently by
computing unbiased estimators of its gradient. Demonstrations with synthetic and
real-world datasets revealed that our approach leads to the best average and
worst-case regret compared to competing methods.
\begin{table}[]
    \centering
    \begin{tabular}{crcccccc}
    \hline \hline\\
     &  & \multicolumn{6}{c}{ $\text{log}_{10}$[normalized simple regret]} \\
    \hline
     Test Case & $d$ & \small{\texttt{EImumax}} & \small{\texttt{PImumax}} & \small{\texttt{UCB}} & \small{\texttt{Random}} & \small{\texttt{R-EI}} & \small{\texttt{r2LEY}} \\
    \hline
    \tiny{Quad-b} & \tiny{$1$}  & \tiny{$-0.32 \pm 0.04$} & \tiny{$-0.47 \pm 0.03$} & \tiny{$-0.50 \pm 0.17$} & \tiny{$-0.31 \pm 0.08$} & \tiny{$-0.48 \pm 0.09$} & \small{$-0.55 \pm 0.06$} \\
    \tiny{Quad-c} & \tiny{$1$}  & \tiny{$-0.57 \pm 0.03$} & \tiny{$-0.93 \pm 0.14$} & \tiny{$-0.89 \pm 0.16$} & \tiny{$-1.26 \pm 0.24$} & \tiny{$-0.90 \pm 0.08$} & \small{$-0.94 \pm 0.07$}\\
    \tiny{Quad-d} & \tiny{$1$}  & \tiny{$-1.45 \pm 0.25$} & \tiny{$-2.06 \pm 0.16$} & \tiny{$-1.29 \pm 0.16$} & \tiny{$-1.41 \pm 0.27$} & \tiny{$-0.94 \pm 0.08$} & \small{$-3.63 \pm 0.27$} \\
    \tiny{Gri} & \tiny{$2$}     & \tiny{$-0.77 \pm 0.04$} & \tiny{$-0.67 \pm 0.04$} & \tiny{$-0.85 \pm 0.042$} & \tiny{$-0.32 \pm 0.04$} & \tiny{$-0.35 \pm 0.04$} & \small{$-1.00 \pm 0.05$} \\
    \tiny{Hart3} & \tiny{$3$}   & \tiny{$-0.37 \pm 0.08$} & \tiny{$-0.33 \pm 0.07$} & \tiny{$-0.46 \pm 0.13$} & \tiny{$-0.51 \pm 0.11$} & \tiny{$-0.75 \pm 0.12$} & \small{$-3.59 \pm 0.39$} \\
    \tiny{Hart6} & \tiny{$6$}   & \tiny{$-0.41 \pm 0.06$} & \tiny{$-0.37 \pm 0.02$} & \tiny{$-0.29 \pm 0.04$} & \tiny{$-0.34 \pm 0.03$} & \tiny{$-0.52 \pm 0.03$} & \small{$-3.31 \pm 0.46$} \\
    \tiny{Levy8} & \tiny{$8$}   & \tiny{$-0.35 \pm 0.02$} & \tiny{$-0.35 \pm 0.02$} & \tiny{$-0.48 \pm 0.04$} & \tiny{$-0.34 \pm 0.01$} & \tiny{$-0.32 \pm 0.02$} & \small{$-1.21 \pm 0.59$} \\
    \tiny{St-Ta} & \tiny{$10$}  & \tiny{$-0.33 \pm 0.02$} & \tiny{$-0.32 \pm 0.02$} & \tiny{$-0.41 \pm 0.04$} & \tiny{$-0.34 \pm 0.02$} & \tiny{$-0.32 \pm 0.01$} & \small{$-0.65 \pm 0.02$} \\
    \tiny{Intel} & \tiny{$2$}   & \tiny{$-0.36 \pm 0.04$} & \tiny{$-0.27 \pm 0.04$} & \tiny{$-0.58 \pm 0.05$} & \tiny{$-0.45 \pm 0.05$} & \tiny{$-0.31 \pm 0.05$} & \small{$-1.24 \pm 0.35$}\\
     \tiny{SARCOS} & \tiny{$7$}  & \tiny{$-0.40 \pm 0.04$} & \tiny{$-0.35 \pm 0.06$} & \tiny{$-0.31 \pm 0.03$} & \tiny{$-0.40 \pm 0.05$} & \tiny{$-0.48 \pm 0.07$} & \small{$-2.52 \pm 0.34$}  \\
    \hline \hline
    
    \end{tabular}
    \caption{Mean $\pm$ standard error of $\T{log}_{10}[\T{normalized simple
    regret}]$ at $T$. Larger fonts represent best method in terms of either
    mean regret or worst-case regret. Each algorithm is repeated 20
    times with independent starting samples to compute metric.}
    \label{tab:results}
\end{table}

\clearpage
\section*{Acknowledgements}
This work was supported by the U.S. Department of Energy, Office of Science, Advanced Scientific Computing Research, Quantum Algorithm Teams Program under Contract DE-AC02-06CH11357.

\bibliography{refs}

\clearpage
\section*{Supplementary material}
\subsection*{Time-Dependent Acquisition Functions}
We present the following time-dependent extensions of existing acquisition functions namely, the EI, PI and UCB.

\subsubsection*{Expected Improvement}
The improvement function is defined as
\begin{equation}
  I(\x,t) = \max \left(0, Y(\x,t)-\xi\right),
    \label{e:improvement}
\end{equation}
where $\xi$ is the target. Then $I(\x, t) \sim \mcl{N}\left((\mu (x, t) - \xi)^+, \sigma^2(\x, t) \right)$. So
\begin{equation}
    \alpha_\T{EI}(\x, t) := \mbb{E}_Y(I) = \int_{I=0}^{\infty} I \left\lbrace \f{1}{\sqrt{2\pi} \sigma(\x, t)} \text{exp} \left[- \f{(I - Y(\x, t) + \xi)^2}{2 \sigma^2(\x, t)} \right]dI \right\rbrace
    \label{e:exp_improvement}
\end{equation}
where the lower limit in the above integral is zero due to the fact that the
improvement is non-negative. Further simplification can be obtained via
transformation of variable as 

\[ u = \f{I-Y+\xi}{\sigma}\] which transforms the integral in \eqref{e:exp_improvement} as
\begin{equation}
    \begin{split}
    \alpha_\T{EI}(\x, t) &= \int_{u=\f{\xi-Y}{\sigma}}^{\infty} (\sigma u + Y - \xi) 
    \f{1}{\sqrt{2\pi}} \text{exp}(-u^2/2) du\\ 
    &= \f{1}{\sqrt{2\pi}} \int_{u=\f{\xi-Y}{\sigma}}^{\infty} \sigma u~ \text{exp}(-u^2/2) du + \f{1}{\sqrt{2\pi}}( Y - \xi) \int_{u=\f{\xi-Y}{\sigma}}^{\infty} \text{exp}(-u^2/2) du \\
    &= \f{1}{\sqrt{2\pi}} \left[\sigma(-\text{exp}(-u^2/2)) \right]_{\f{\xi-Y}{\sigma}}^{\infty} + \f{1}{\sqrt{2\pi}} ( Y - \xi) \int_{u=\f{\xi-Y}{\sigma}}^{\infty} \text{exp}(-u^2/2) du
    \end{split}
    \label{e:exp_improvement2}
\end{equation}
By denoting the standard normal pdf and cdf as $\phi()$ and $\Phi()$
respectively, the above equation can be re-written concisely as 

\begin{equation}
   \alpha_\T{EI}(\x, t) = \sigma(\x, t) \Phi\left(\f{\xi-\mu}{\sigma}\right) + ( \mu - \xi) \phi \left( \f{\xi-\mu}{\sigma}\right)
\end{equation}

Note that for time-dependent problems, where the oracle is changing with time, it is required that $\xi = \xi(t)$. 

\subsubsection*{GP Upper Confidence Bound}
The GP-UCB for time-dependent oracles is specified as an optimistic estimate of the posterior GP mean as follows
\begin{equation}
    \alpha_\T{UCB}(\x, t) = \mu(\x, t) + \beta^{1/2} \sigma(\x, t),
\end{equation}
\subsection*{Probability of Improvement}
The PI acquisition function for time-dependent oracles is given as follows
\begin{equation}\label{eq:acquisition_PI}
    \begin{split}
        \alpha_\T{PI}(\x, t) =& P(Y(\x, t) \geq \xi) \\
        =& \Phi \left( \f{\mu(\x, t) - \xi}{\sigma(\x)} \right), 
    \end{split}
\end{equation}
where as in EI, the target $\xi = \xi(t)$ for time-dependent problems.
\section*{1D Quadratic Test Functions}
\begin{figure}[ht]
    \centering
        \begin{subfigure}[htb!]{\textwidth}
        \includegraphics[width=0.4\textwidth, trim= 8cm 1cm 8cm 1cm, clip]{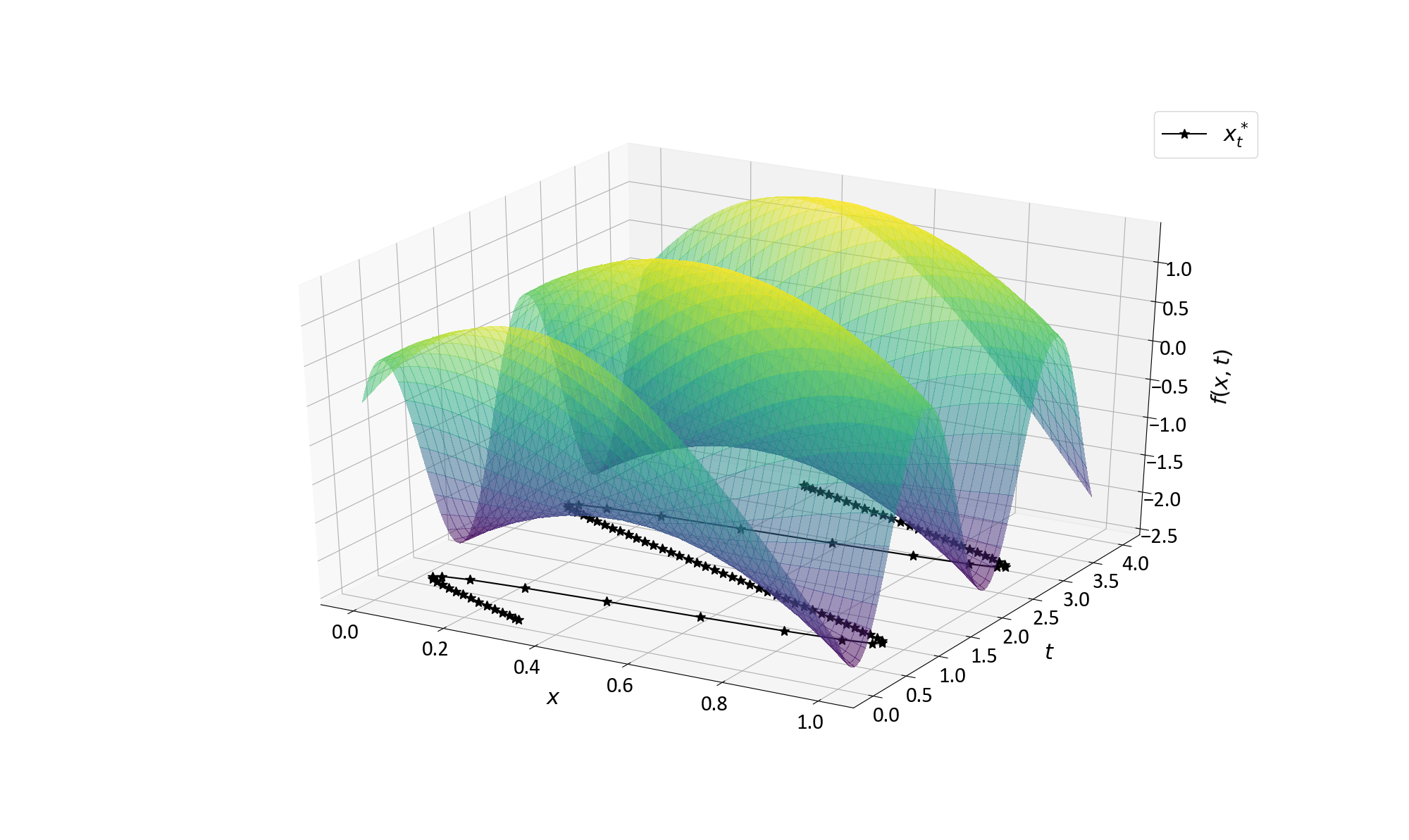}%
        \hfill
        \includegraphics[width=0.4\textwidth]{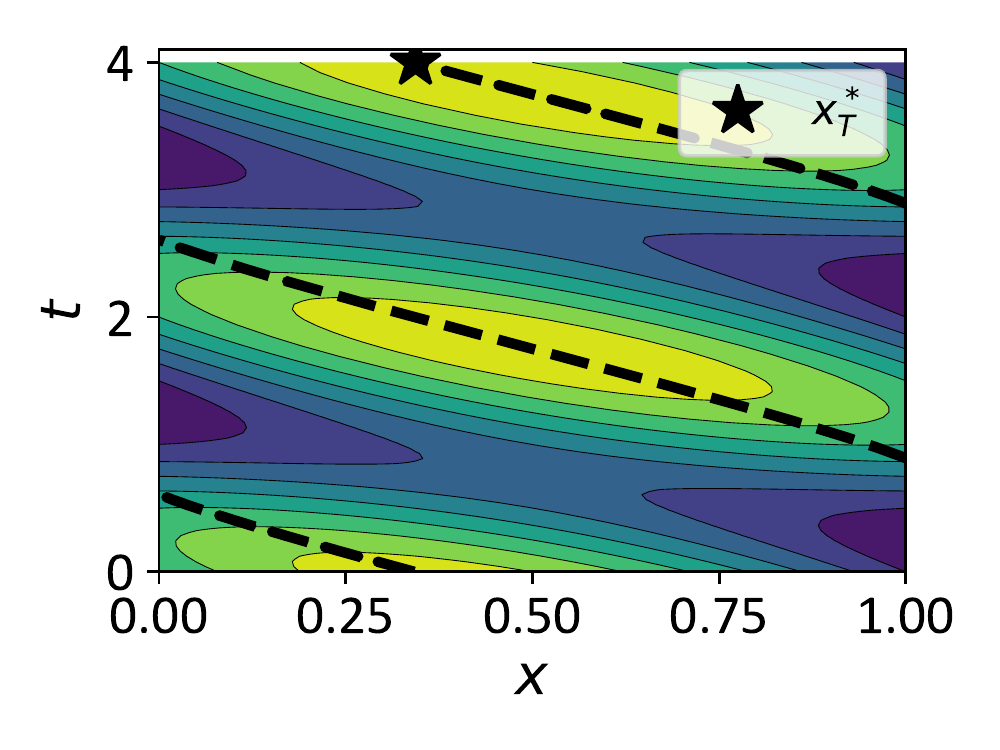}%
        \caption{Quadratic-a}
        \label{sf:}
        \end{subfigure}\\
        \begin{subfigure}[htb!]{\textwidth}
        \includegraphics[width=0.4\textwidth, trim= 8cm 1cm 8cm 1cm, clip]{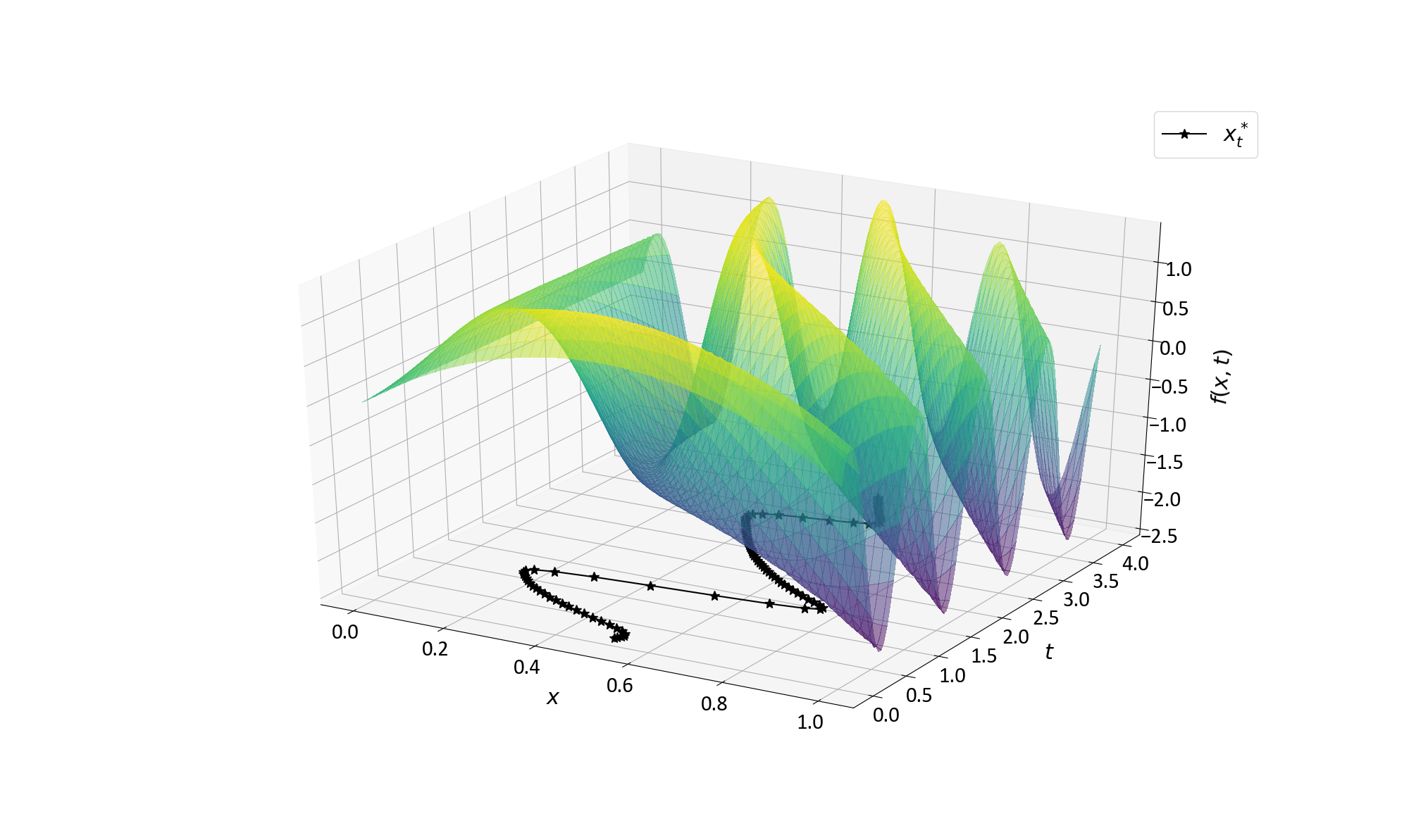}%
        \hfill
        \includegraphics[width=0.4\textwidth]{quadratic_f_xt_b_contour.pdf}%
        \caption{Quadratic-b}
        \label{sf:}
        \end{subfigure}\\
        \begin{subfigure}[htb!]{\textwidth}
        \includegraphics[width=0.4\textwidth, trim= 8cm 1cm 8cm 1cm, clip]{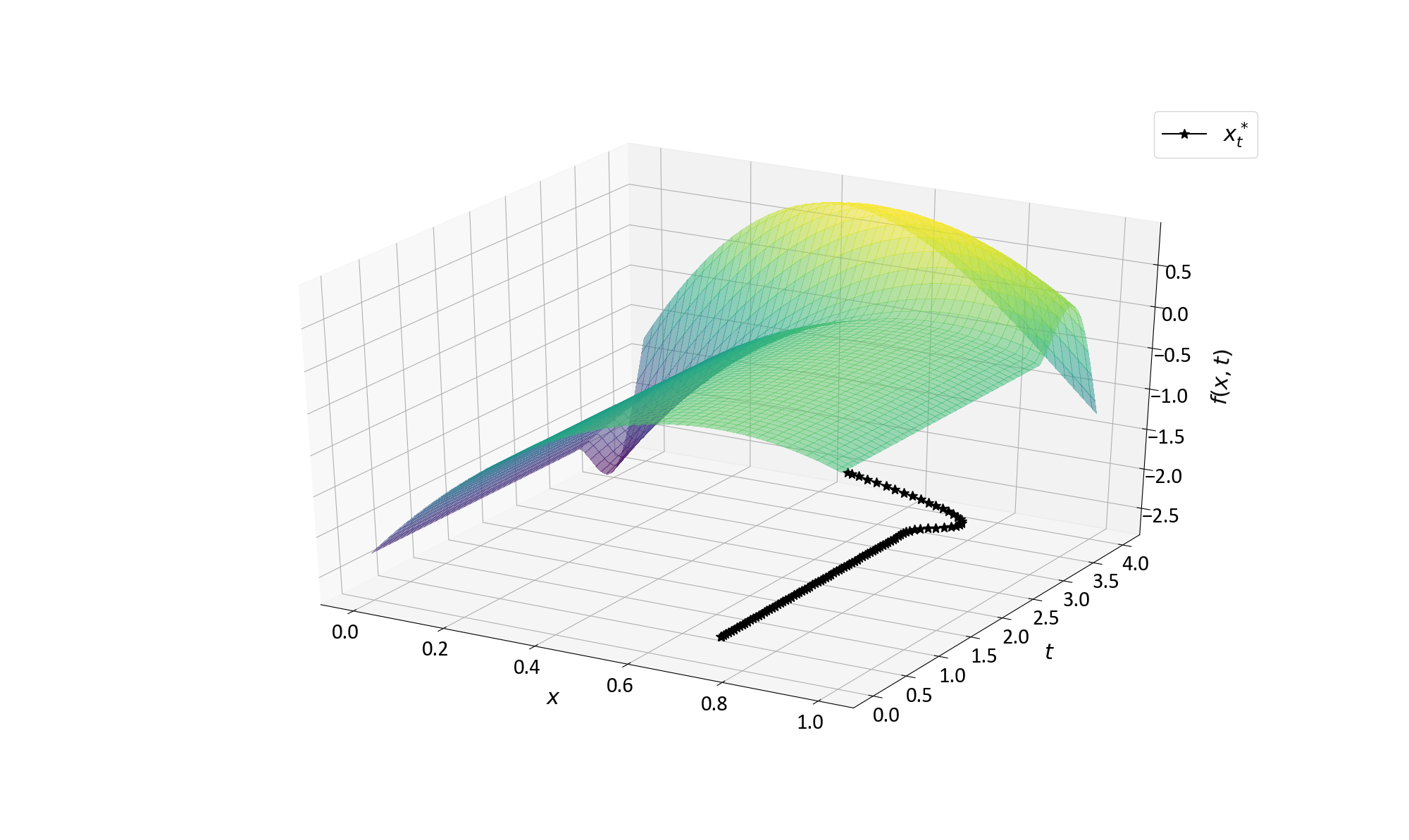}%
        \hfill
        \includegraphics[width=0.4\textwidth]{quadratic_f_xt_c_contour.pdf}%
        \caption{Quadratic-c}
        \label{sf:}
        \end{subfigure}\\
        \begin{subfigure}[htb!]{\textwidth}
        \includegraphics[width=0.4\textwidth, trim= 8cm 1cm 8cm 1cm, clip]{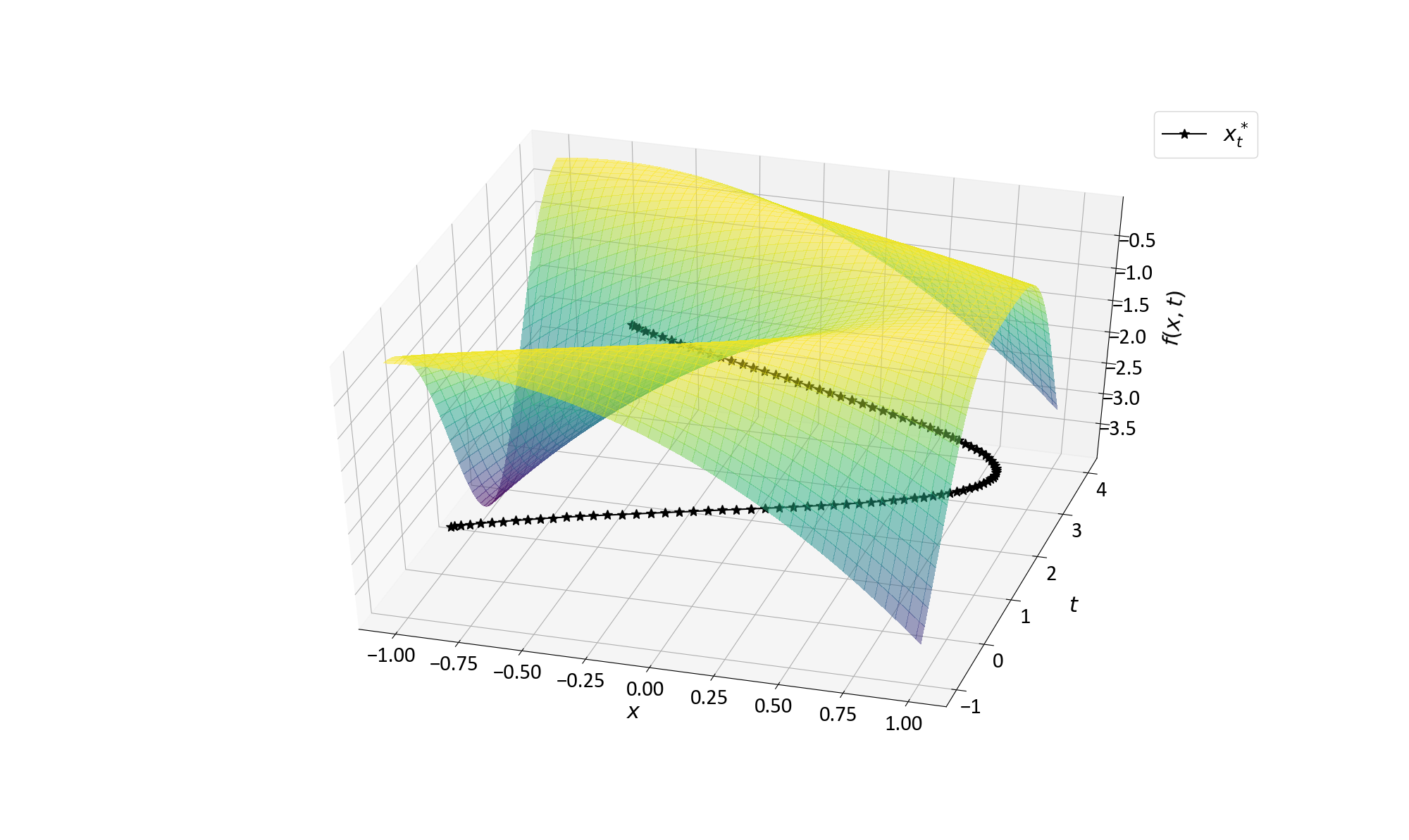}%
        \hfill
        \includegraphics[width=0.4\textwidth]{quadratic_f_xt_d_contour.pdf}%
        \caption{Quadratic-d}
        \label{sf:}
        \end{subfigure}\\        
    \caption{Quadratic ($d=1$) test cases. Dashed lines in the contour plots represent trajectory of local maximizer. Notice that in Qudratic-a, the maximizer makes sudden jumps (from one boundary to other) at around $t=1$ and $t=3$ and, Quadratic-b has two local maxima starting $t=2$.} 
    \label{f:quads}
\end{figure}

We first consider quadratic 1D functions $f$
over the domain $\mcl{X} = [0,1]$ with
\begin{equation} 
f_{\x}(\x) = -\alpha \times (\x - s)^2,
\label{eq:1dq}
\end{equation}
where $s$ and $\alpha$ are scalar parameters. The context-dependent component takes the three forms
\begin{equation*}
    \begin{split}
    \T{Quadratic-a}: \quad f_{\x t} = &\T{sin}\left(\pi(\x+t)\right) + \T{cos}\left(\pi(\x+t)\right) \\
    \T{Quadratic-b}: \quad f_{\x t} = &\T{sin}\left(\pi(\x t)\right) + \T{cos}\left(\pi(\x t)\right) \\
    \T{Quadratic-c}: \quad f_{\x t} = &\T{sin}\left(\pi(\x[t-3]^+)\right) +
    \T{cos}\left(\pi(\x[t-3]^+)\right)\\
    \T{Quadratic-d}: \quad f_{\x t} = &2\x\T{sin}(t) - \T{sin}^2(t),  
    \end{split}
\end{equation*}
where, $[t-3]^+ \equiv \T{max}(0, t-3)$ is specified to induce movement of
$\x_t^*$ for $t\geq 3$. The different $f_{\x t}$ functions are chosen to create
varying patterns of context-dependent $\x_t^*$. In \eqref{eq:1dq}, we set $s=0.5, \alpha=4.0$; see Figure \ref{f:quads} for further details.

\subsection*{Details of synthetic test functions}

\begin{table}[htb!]
    \centering
        \caption{Domain and maximizer of the time-independent part of the synthetic test functions. Note that each test function is composed of a time-dependent component $f_{\x t}$ (not included in the table) that moves the maximizer with respect to $t$ in the experiments}
    \begin{tabular}{ccccc}
    \hline
    Test case & $d$ & $\mcl{X}$ & $\mcl{T}$ & $\underset{\x \in \mcl{X}}{arg\T{max}}~f_{\x}(\x)$ \\
    \hline \hline
    Quadratic-a & $1$ & $[0, 1]$ & $[0, 4]$ & $0.5$ \\
    Quadratic-b & $1$ & $[0, 1]$ & $[0, 4]$ & $0.5$ \\
    Quadratic-c & $1$ & $[0, 1]$ & $[0, 4]$ & $0.5$ \\
    Quadratic-d & $1$ & $[0, 1]$ & $[0, 4]$ & $0.5$ \\
    Griewank & $2$ & $[-5, 5]^2$ & $[0, 4]$ & $[0, 0]$ \\
    Hartmann-3d & $3$ & $[0, 1]^3$ & $[0, 4]$ & $[0.11, 0.56, 0.85]$ \\
    Hartmann-6d & $6$ & $[0, 1]^6$ & $[0, 4]$ & $[0.20,0.15,0.48,0.28,0.31,0.66]$ \\    
    Levy & $8$ & $[-10, 10]^8$ & $[0, 4]$ & $[1,\ldots,1]$ \\  
    Styblinski-Tang & $10$ & $[-5, 5]^{10}$ & $[0, 4]$ & $[-2.903534, \ldots, -2.903534]$ \\  
    Intel Sensor & $2$ & $[0.5, 40.5]^2$ & $[0, 100]$ & - \\  
    SARCOS Robot & $7$ & $[0, 2]^7$ & $[0, 4]$ & - \\  
    
    \hline
    \end{tabular}
    \label{tab:my_label}
\end{table}

\clearpage

\subsection*{Proof of Theorem 3.1}
\begin{theorem}[Interchange of gradient and expectation operators]  Let us write
  $g^*(\x_{j+1}, T)|y_{j+1}$ as $g^*(\x_{j+1}, T, y_{j+1})$ and, let $\mcl{Y}$
  be the support of $y_{j+1}$, whose density is $p(y_{j+1})$, and the
  domain of $\x_{j+1}$ be $\mcl{X}$. Let $g^*(\x_{j+1}, T, y_{j+1})\times
  p(y_{j+1})$ and $\partial g^*(\x_{j+1}, T, y_{j+1})/\partial \x_{j+1} \times
  p(y_{j+1})$ be continuous on $\mcl{X} \times \mcl{Y}$. Furthermore, assume that the kernel $k()$ is continuously differentiable in $\mcl{X}$. Suppose that there
  exist nonnegative functions $q_0(y_{j+1})$ and $q_1(y_{j+1})$ such that
  $|g^*(\x_{j+1}, T, y_{n+1}) \times p(y_{j+1})| \leq q_0(y_{j+1})$ and $\|
  \partial g^* / \partial \x_{j+1} \times  p(y_{j+1})\| \leq q_1(y_{j+1})$ for
  all $(\x_{j+1}, y_{j+1}) \in \mcl{X} \times \mcl{Y}$, where
  $\int_{\mcl{Y}}q_1(y_{j+1}) dy_{j+1} < \infty$ and
  $\int_{\mcl{Y}}q_2(y_{j+1}) dy_{j+1} < \infty$. Then,
\[ \nabla \mbb{E}_j \left[ g^*(\x_{j+1}, T, y_{j+1}) \right] = \mbb{E}_j \left[\nabla g^*(\x_{j+1}, T, y_{j+1}) \right]\]

Since $\nabla \alpha_{2LEY}(\x_{j+1}) = \mbb{E}_j \left[\nabla g^*(\x_{j+1}, T,
y_{j+1}) \right]$, a realization of $\nabla g^*(\x_{j+1}, T, y_{j+1})$ yields
an unbiased estimate of the true gradient.
\end{theorem}
In what follows, we use $g^*(\x_{n+1})$ to denote  $g^*(\x_{n+1}, T, y_{n+1})$ for the sake of brevity
\begin{equation}
\begin{split}
    \frac{\partial}{\partial \x_{n+1}} \mbb{E}_n \left[ g^*(\x_{n+1}) \right] &= \underset{h\rightarrow 0}{\T{lim}}~\frac{1}{h} \left \lbrace \mbb{E}_n \left[ g^*(\x_{n+1}+h) \right] -  \mbb{E}_n \left[ g^*(\x_{n+1}) \right] \right \rbrace \\
    &= \underset{h\rightarrow 0}{\T{lim}}~ \left \lbrace \mbb{E}_n \frac{1}{h} \left[ g^*(\x_{n+1}+h) - g^*(\x_{n+1}) \right] \right \rbrace \\
    &= \underset{h\rightarrow 0}{\T{lim}}~ \left \lbrace \mbb{E}_n \frac{\partial g^*(\bar{\x}_{n+1})}{\partial \x_{n+1}}  \right \rbrace ,\\
\end{split}    
\end{equation}
where, $\bar{\x}_{n+1} = \lambda \x_{n+1} + (1-\lambda)(\x_{n+1} + h)$ for some $\lambda \in [0,1]$ and the last line follows from the first order mean-value theorem by assuming (but proven below) that $\nabla g^*(\x_{n+1})$ exists and is continuous in $\mcl{X}$.

Finally, we bring the limit inside the integral by the Lebesgue's dominated convergence theorem and state
\begin{equation}
    \frac{\partial}{\partial \x_{n+1}} \mbb{E}_n \left[ g^*(\x_{n+1}) \right] = \left \lbrace \mbb{E}_n~\underset{h\rightarrow 0}{\T{lim}} \frac{\partial g^*(\bar{\x}_{n+1})}{\partial \x_{n+1}}  \right \rbrace = \mbb{E}_n \frac{\partial g^*(\x_{n+1})}{\partial \x_{n+1}}.
    \label{e:LDC}
\end{equation}

To prove that $g^*(\x_{j+1}, T, y_{j+1})\times
  p(y_{j+1})$ and $\partial g^*(\x_{j+1}, T, y_{j+1})/\partial \x_{j+1} \times
  p(y_{j+1})$ 
are continuous on $\mcl{X} \times \mcl{Y}$, first note that $y_{n+1}\sim \mcl{N}(\mu_{n+1}, \sigma^2_{n+1})$ and hence its density $p(y)$ is continuous in $\mcl{Y}$ owing to its exponential structure. 

Consider the expansion of $g^*$ and $\nabla g^*$ repeated below
\begin{equation}
    g^*(\x_{n+1}, T) = \mbf{k}_{n+1}^\top \mbf{K}^{-1}_{n+1} \mbf{y}_{n+1},
    \label{e:2LEY}
\end{equation}
\begin{equation}
        \nabla g^*(\x_{n+1}, T) = \f{\partial \mbf{k}_{n+1}^\top}{\partial \x_{n+1}} \mbf{K}^{-1}_{n+1} +
        \mbf{k}_{n+1}^\top~\mbf{K}^{-1}_{n+1} \f{\partial \mbf{K}_{n+1}}{\partial \x_{n+1}}\mbf{K}^{-1}_{n+1}.
        \label{e:2LEY-grad}
\end{equation}
Since each element of $\mbf{k}_{n+1}$ is the output of a continuously differentiable kernel (by assumption), $\mbf{k}_{n+1}$ is continuously differentiable in $\mcl{X}$. Similarly, due to the differentiability of the kernel, $\mbf{K}_{n+1}$ is continuously differentiable and so is its inverse. Since the product of continuous functions is continuous, $g^*(\x_{n+1})$ is continuous in $\mcl{X}$. By the same arguments and by the fact that the elementwise derivative in $\f{\partial \mbf{K}_{n+1}}{\partial \x_{n+1}}$ is continuous, $\nabla g^*$ is also continuous in $\mcl{X}$.

Finally, since products of continuous functions are continuous, the products $g^*(\x_{n+1}, y_{n+1}) p_{\mcl{Y}}(y_{n+1})$ and $\partial g^*(\x_{n+1}, y_{n+1})/\partial \x_{n+1} p_{\mcl{Y}}(y_{n+1})$ are continuous on $\mcl{X} \times \mcl{Y}$. Therefore \eqref{e:LDC} holds. $\square$

\subsection*{Experiments}
\begin{figure}[ht]
    \centering
        \begin{subfigure}[htb!]{\textwidth}
        \includegraphics[width=0.5\textwidth]{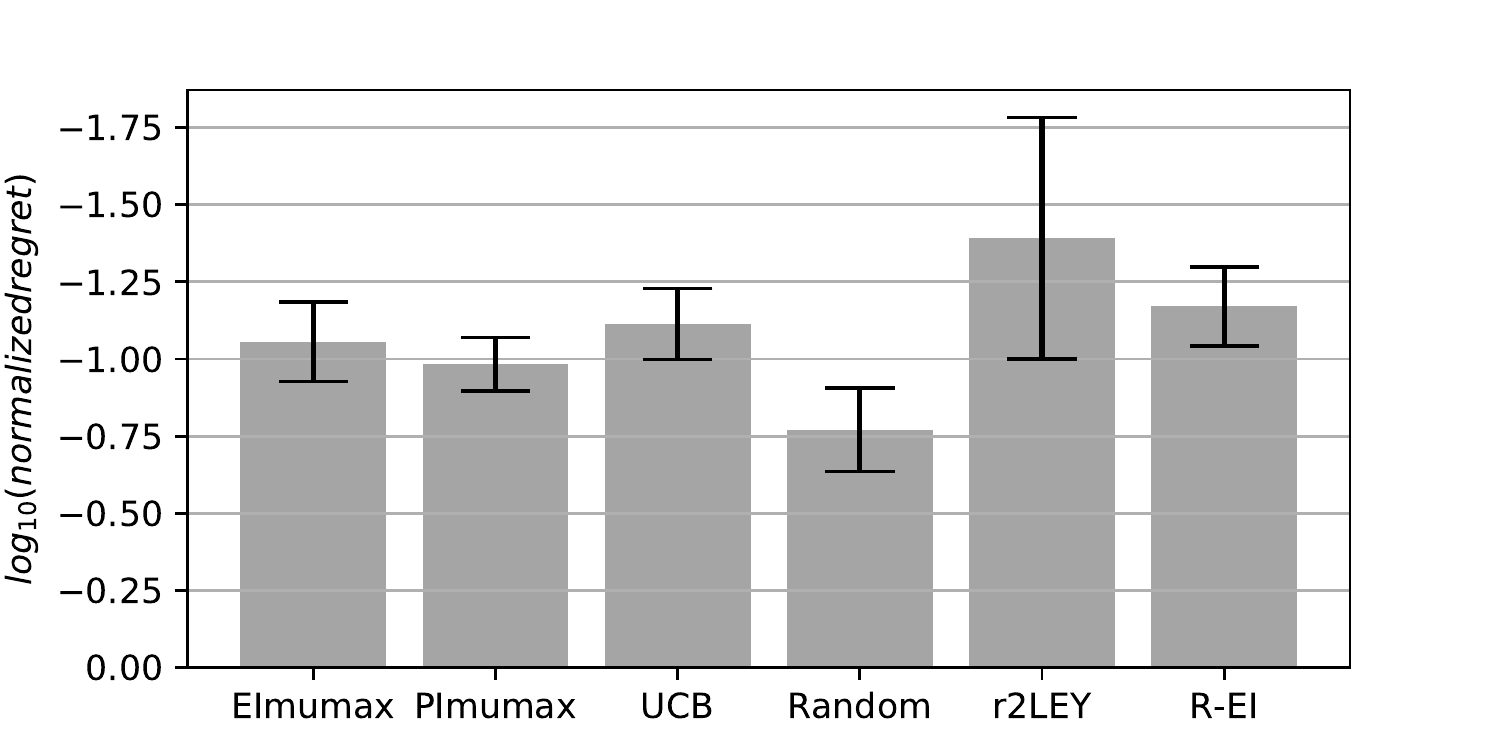}%
        \hfill
        \includegraphics[width=0.5\textwidth]{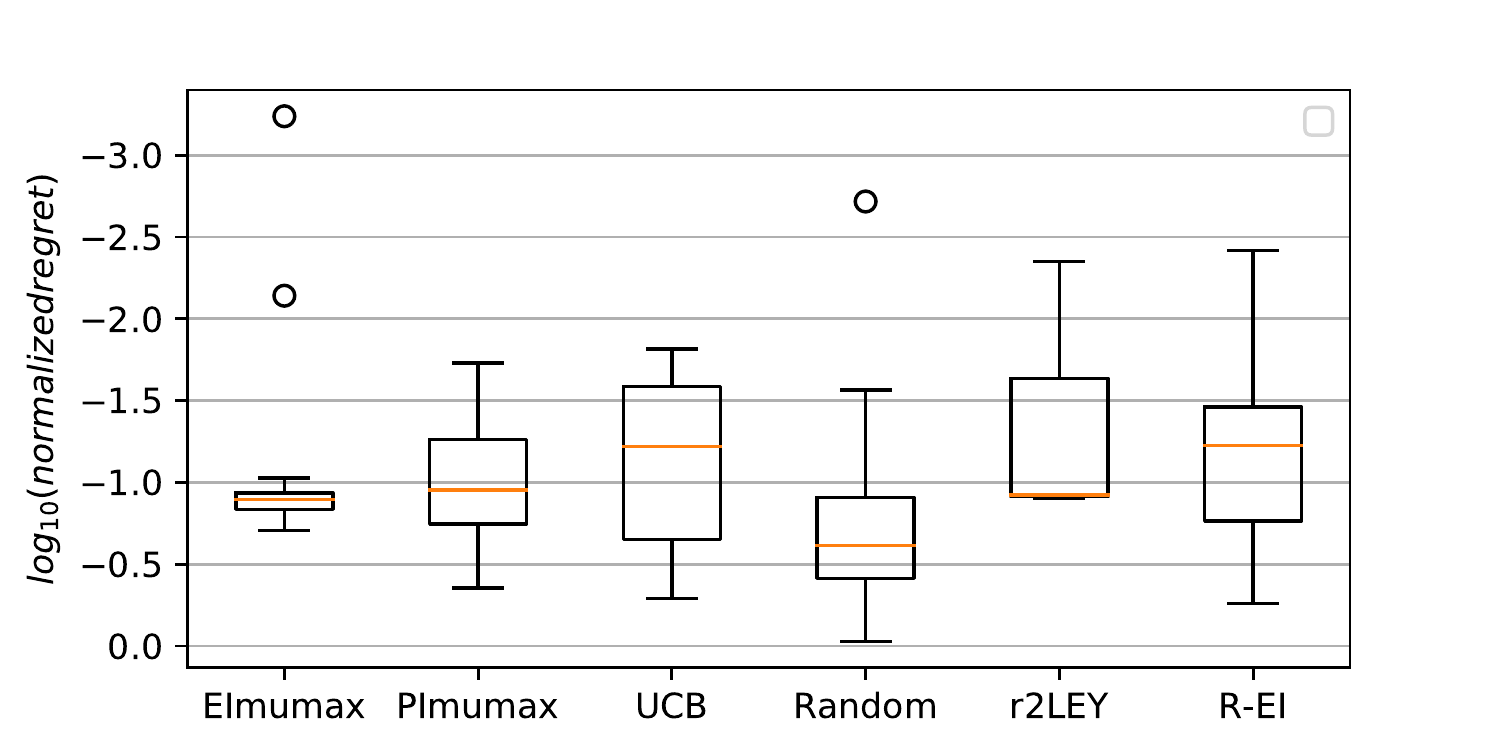}%
        \caption{Quadratic-a}
        \label{sf:}
        \end{subfigure}\\
        \begin{subfigure}[htb!]{\textwidth}
        \includegraphics[width=0.5\textwidth]{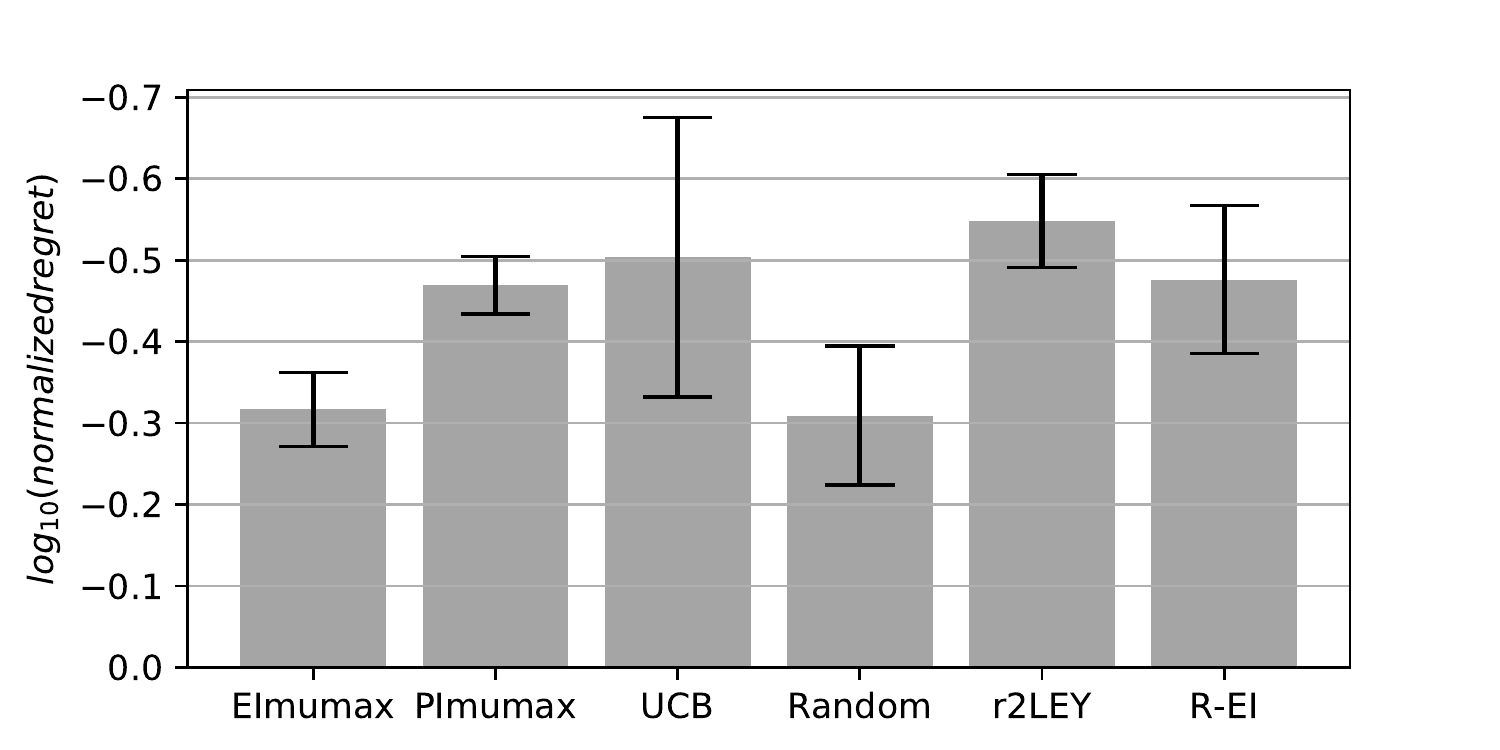}%
        \hfill
        \includegraphics[width=0.5\textwidth]{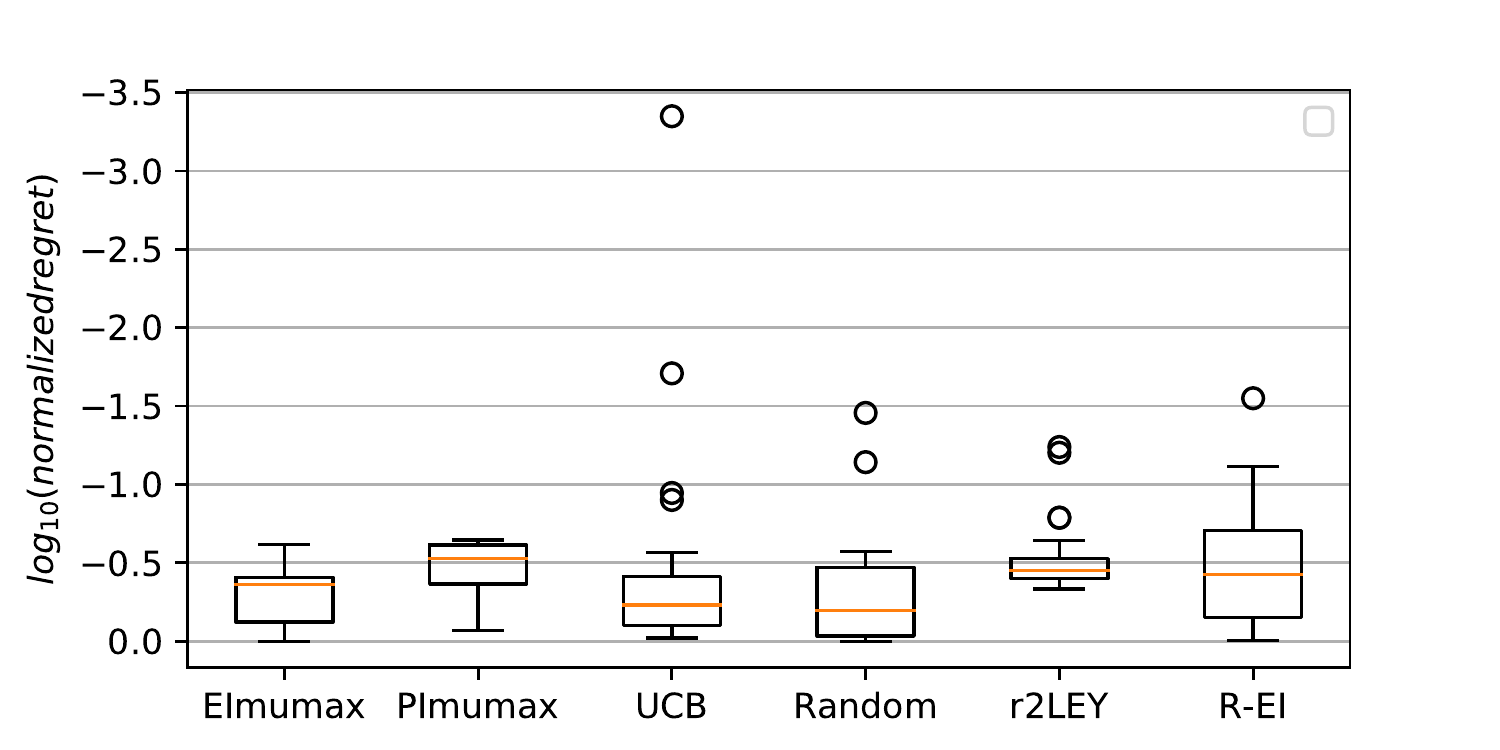}%
        \caption{Quadratic-b}
        \label{sf:}
        \end{subfigure}\\
        \begin{subfigure}[htb!]{\textwidth}
        \includegraphics[width=0.5\textwidth]{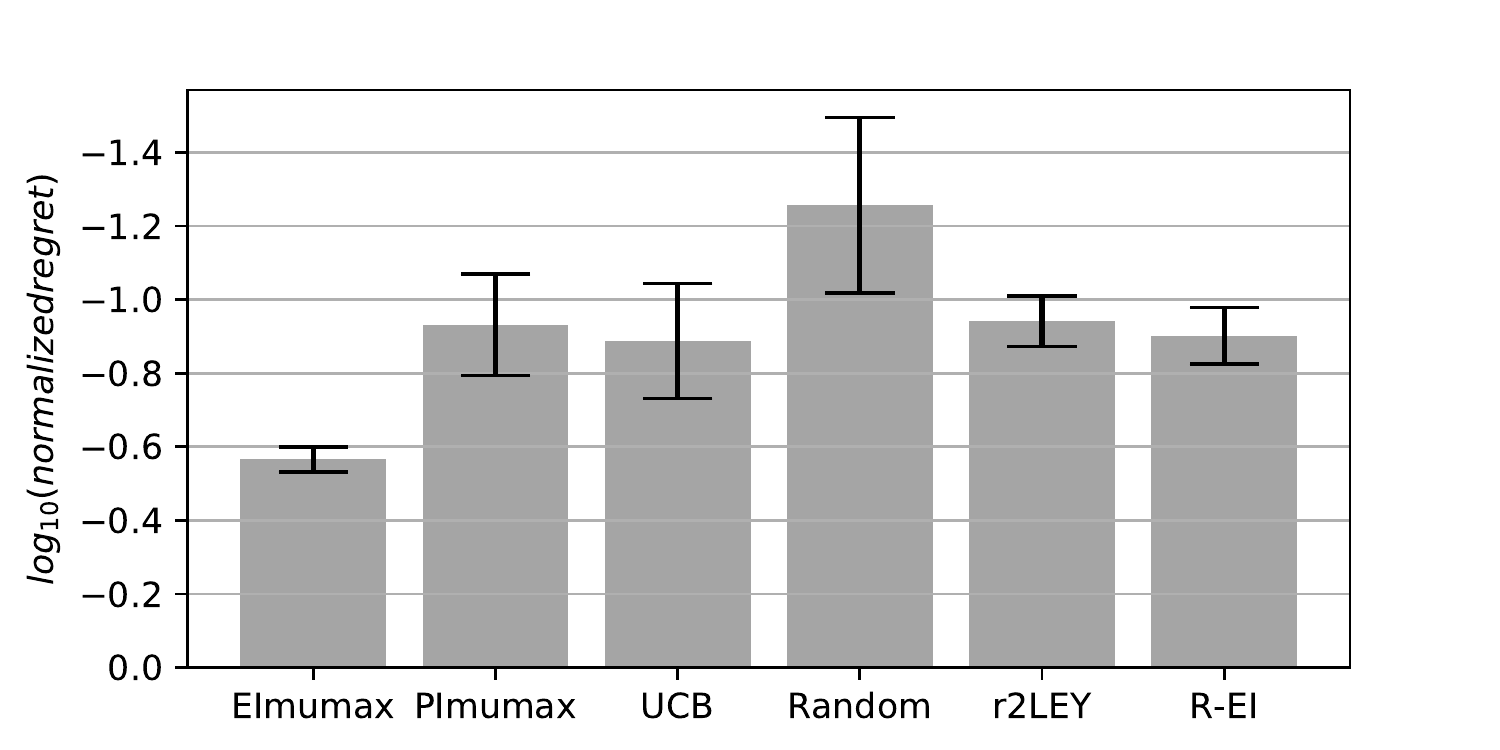}%
        \hfill
        \includegraphics[width=0.5\textwidth]{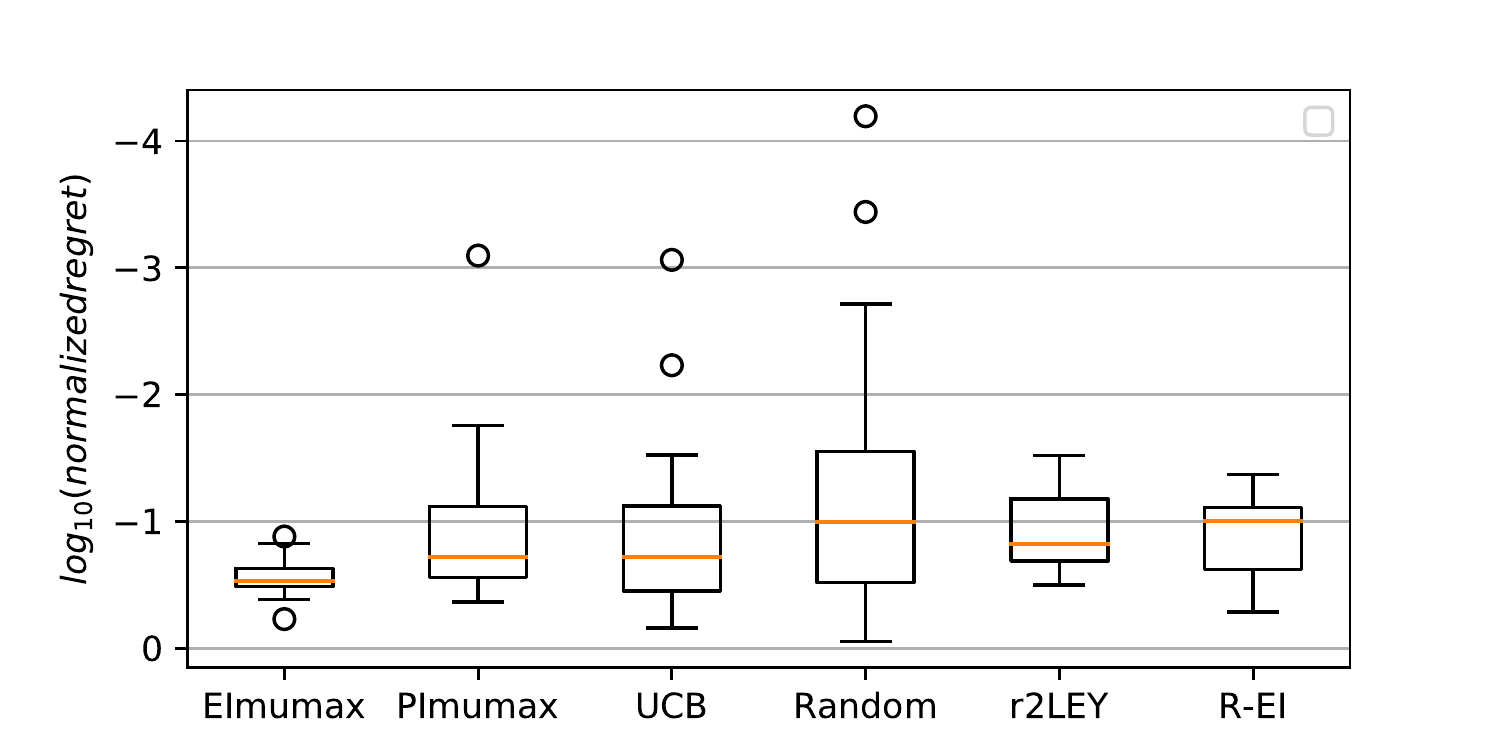}%
        \caption{Quadratic-c}
        \label{sf:}
        \end{subfigure}\\        
        \begin{subfigure}[htb!]{\textwidth}
        \includegraphics[width=0.5\textwidth]{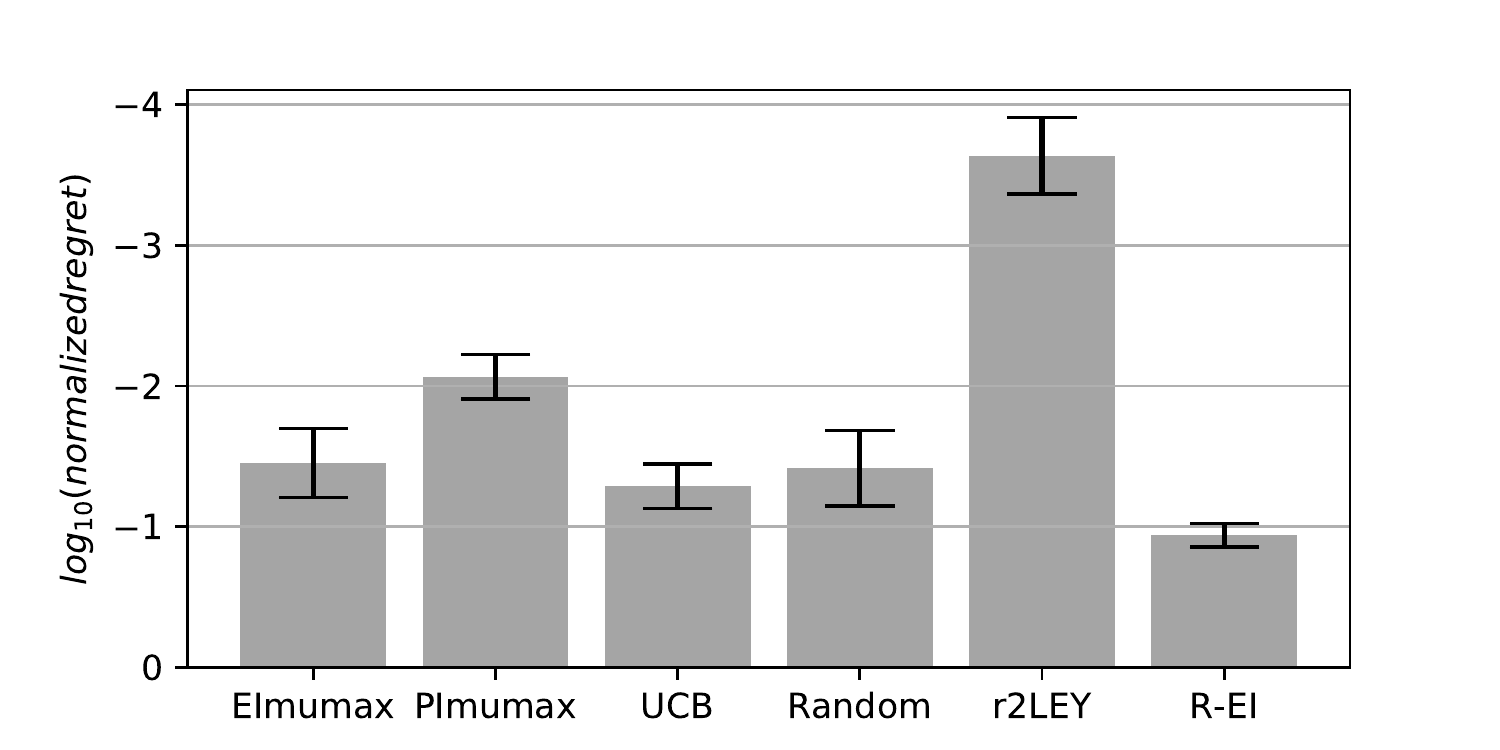}%
        \hfill
        \includegraphics[width=0.5\textwidth]{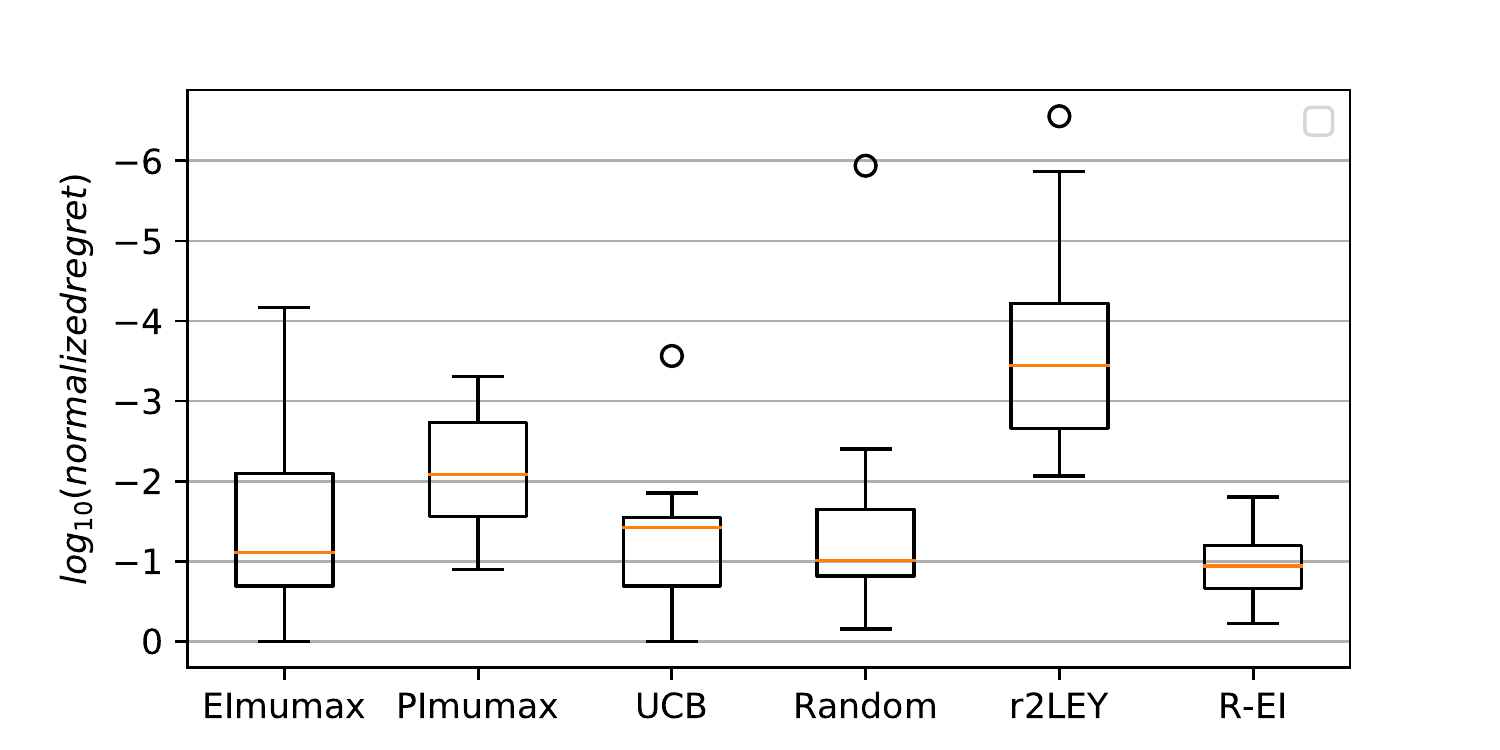}%
        \caption{Quadratic-d}
        \label{sf:}
        \end{subfigure}\\
        \begin{subfigure}[htb!]{\textwidth}
        \includegraphics[width=0.5\textwidth]{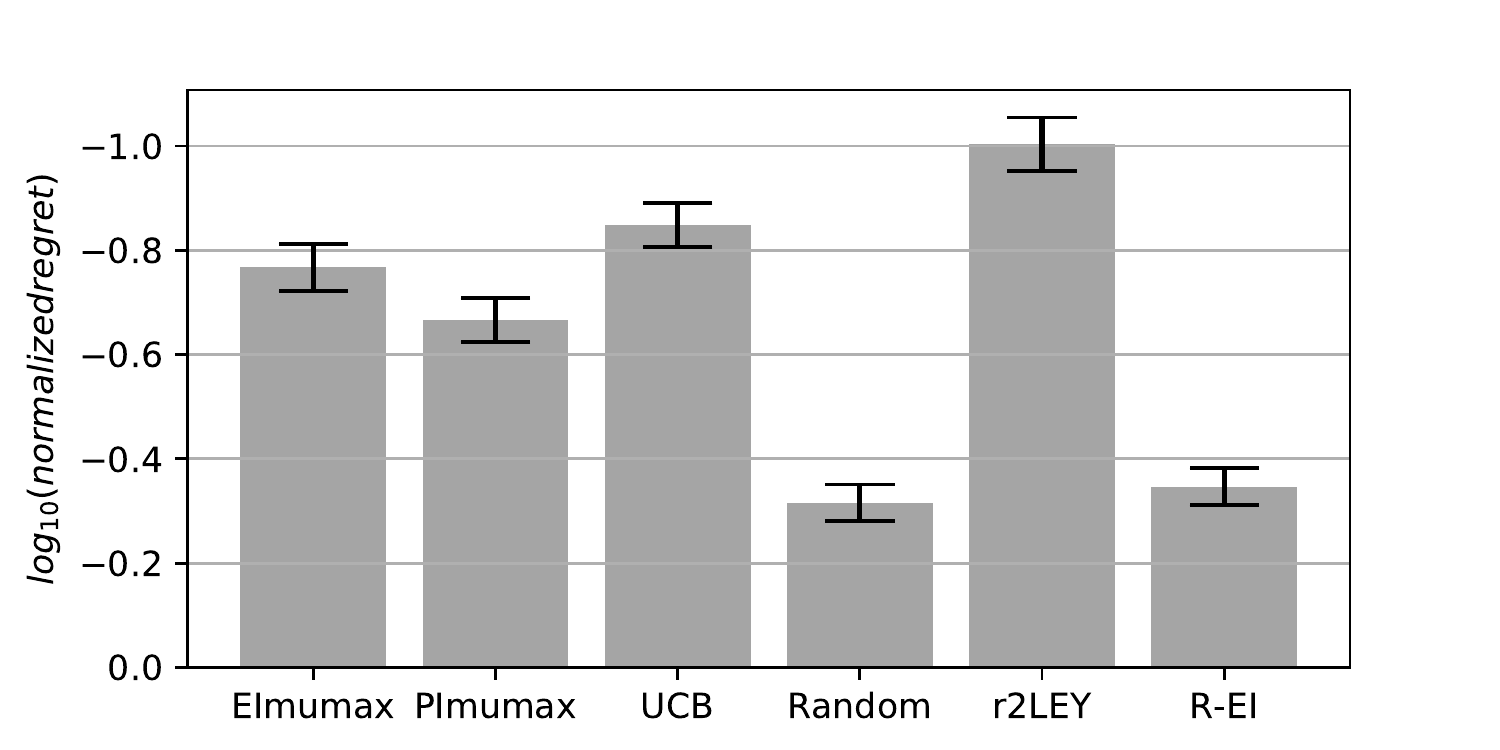}%
        \hfill
        \includegraphics[width=0.5\textwidth]{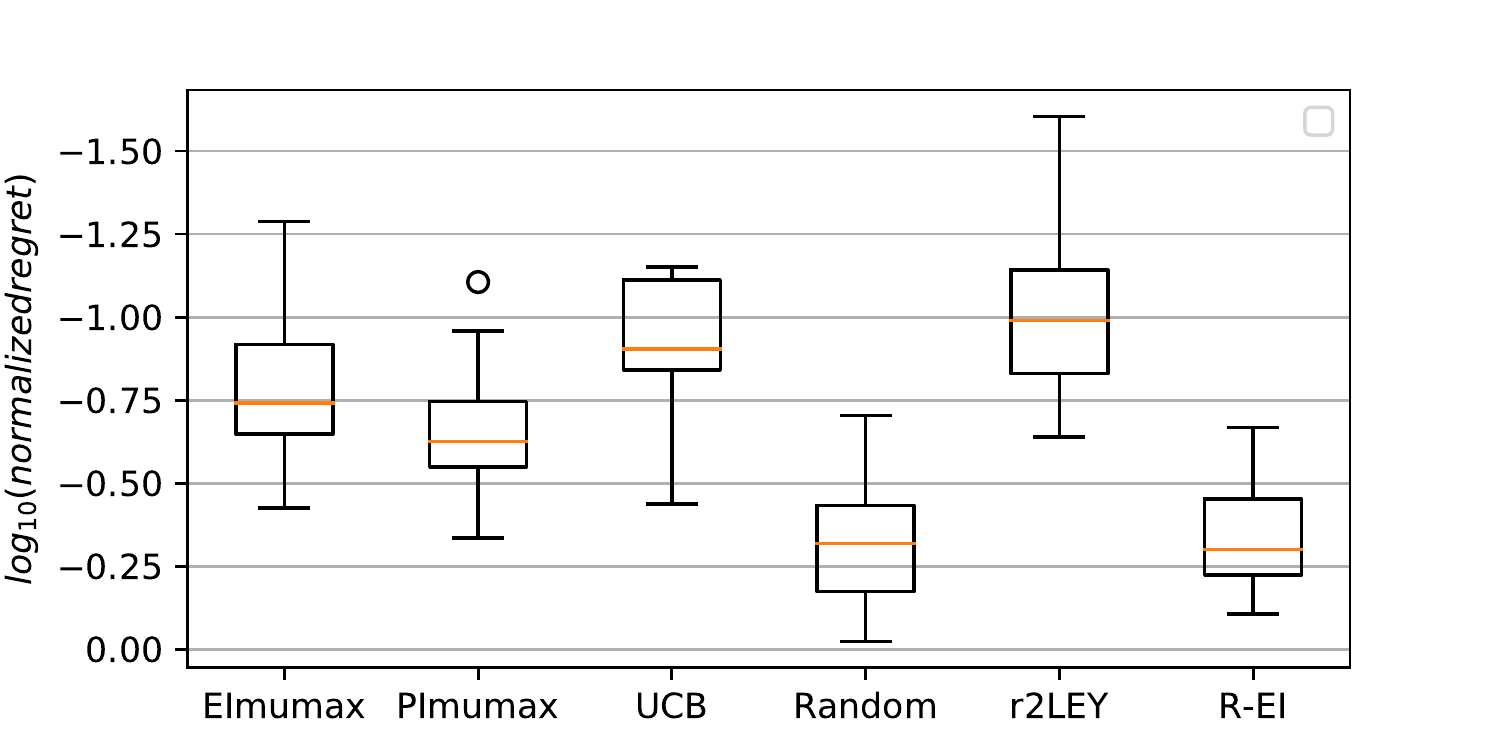}%
        \caption{Griewank-2d}
        \label{sf:}
        \end{subfigure}\\
        
    \caption{} 
    \label{f:syn_results}
\end{figure}

\begin{figure}[ht]
    \ContinuedFloat
    \centering
        \begin{subfigure}[htb!]{\textwidth}
        \includegraphics[width=0.5\textwidth]{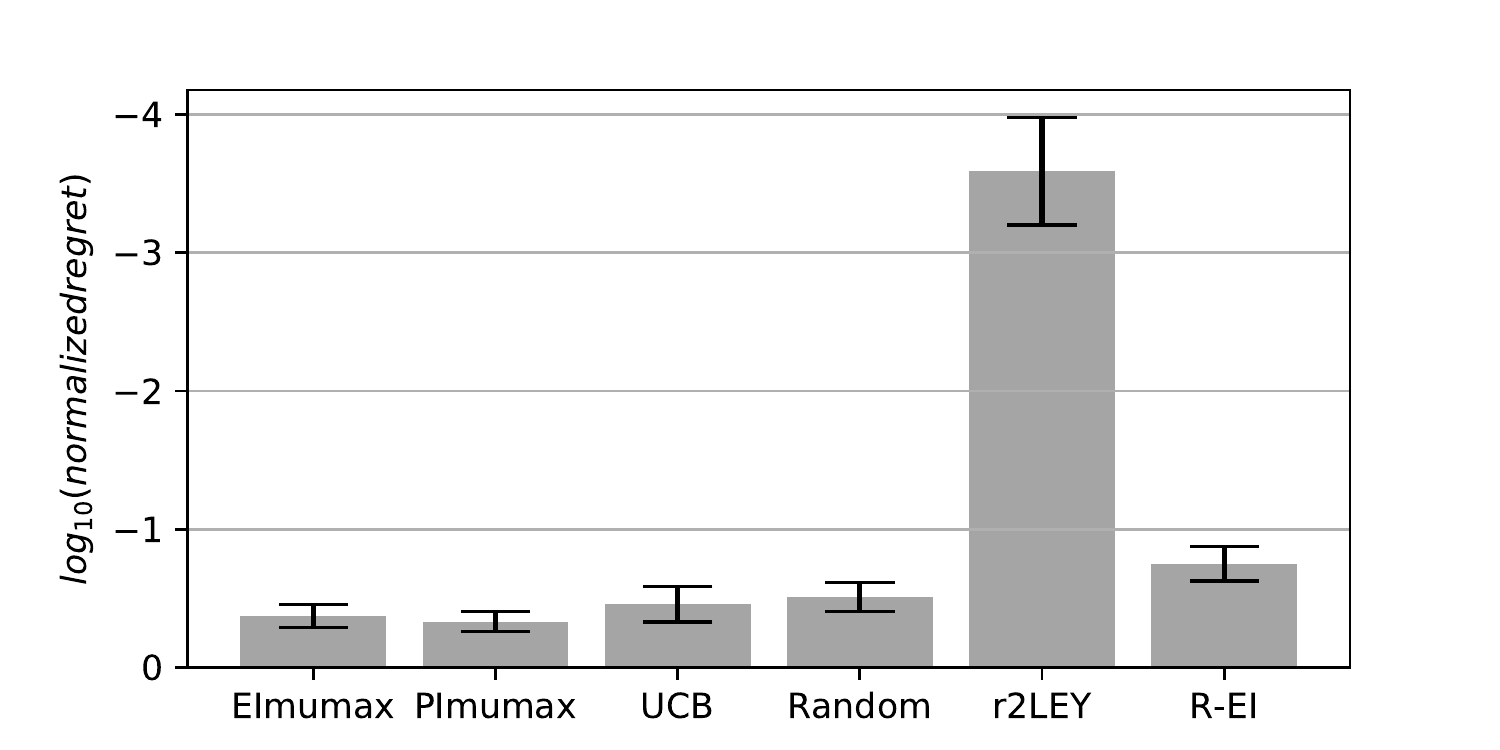}%
        \hfill
        \includegraphics[width=0.5\textwidth]{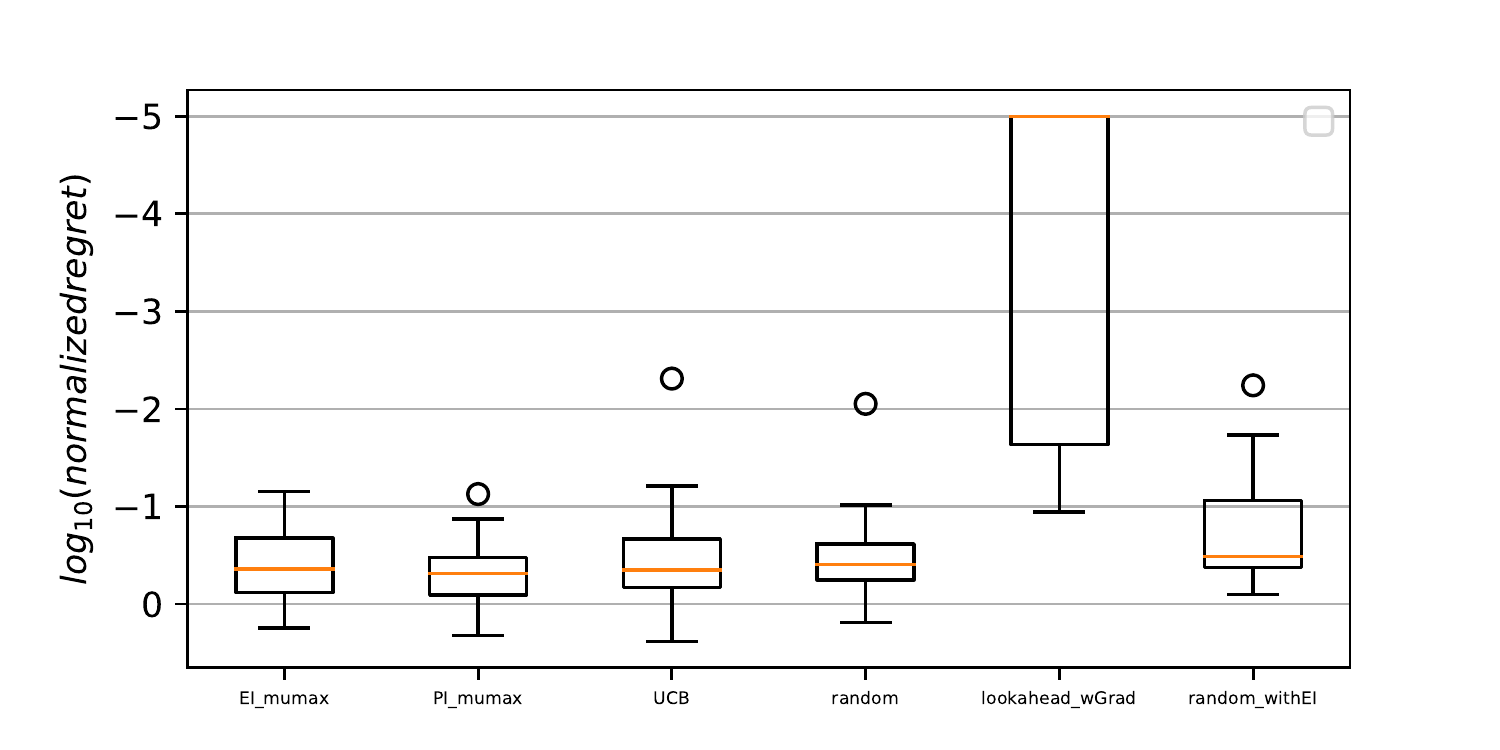}%
        \caption{Hartmann-3d}
        \label{sf:}
        \end{subfigure}\\
        \begin{subfigure}[htb!]{\textwidth}
        \includegraphics[width=0.5\textwidth]{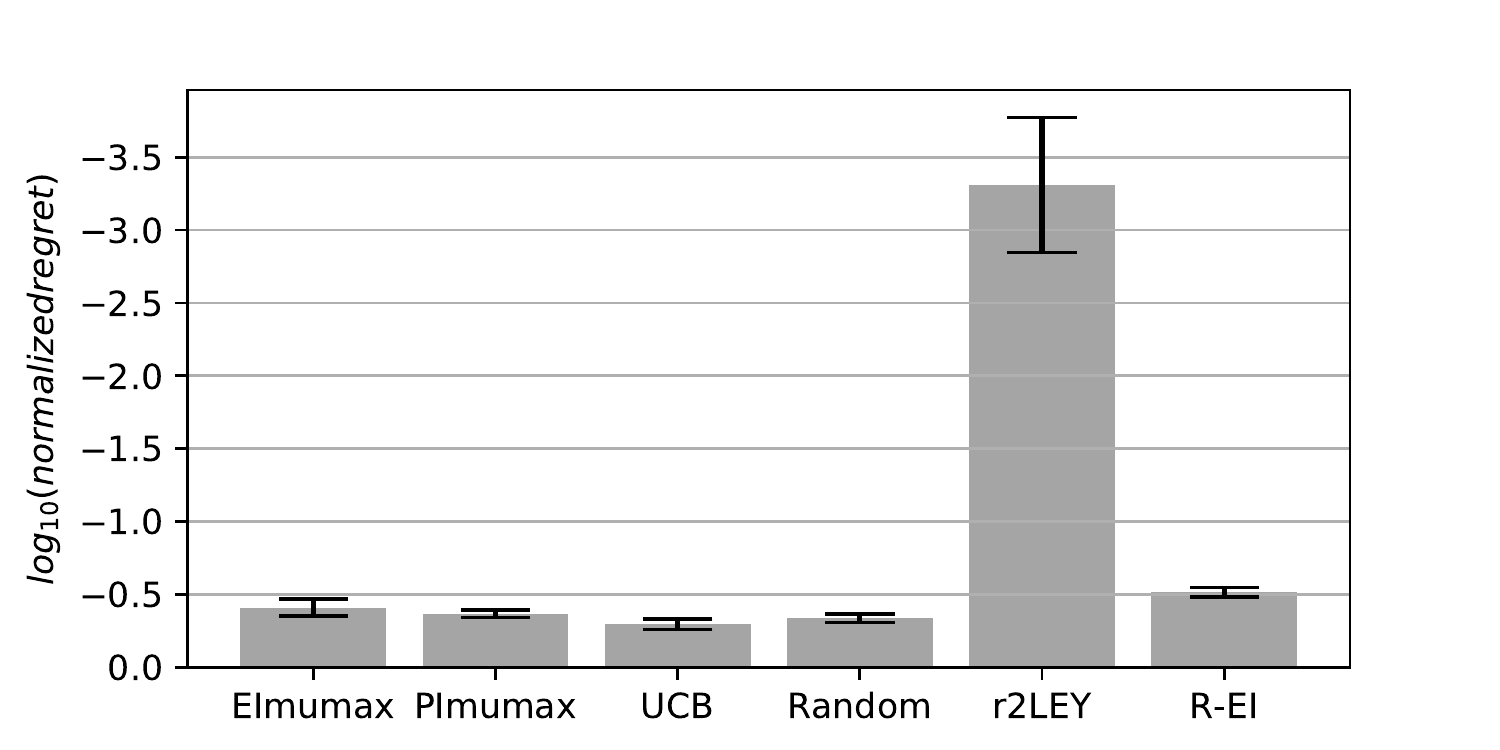}%
        \hfill
        \includegraphics[width=0.5\textwidth]{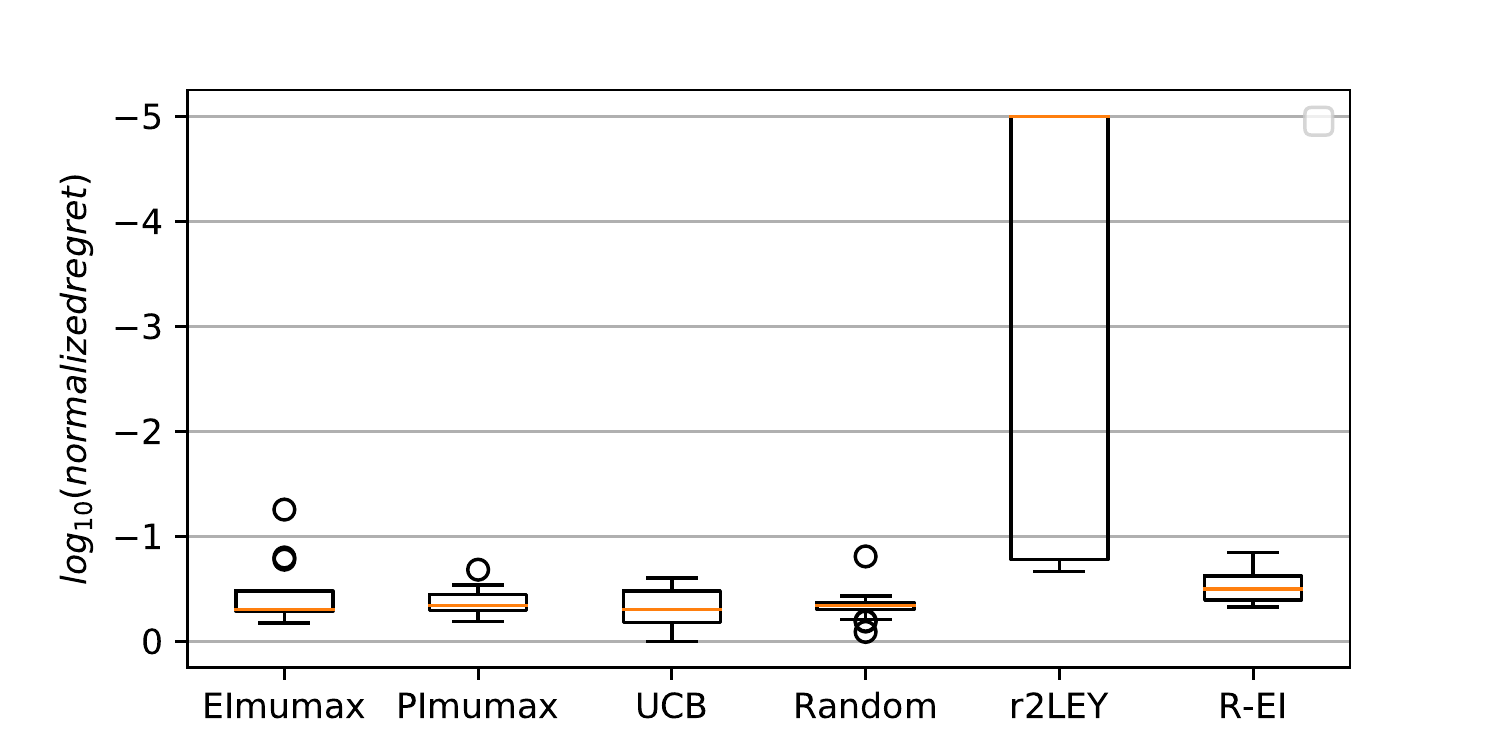}%
        \caption{Hartmann-6d}
        \label{sf:}
        \end{subfigure}\\        
        \begin{subfigure}[htb!]{\textwidth}
        \includegraphics[width=0.5\textwidth]{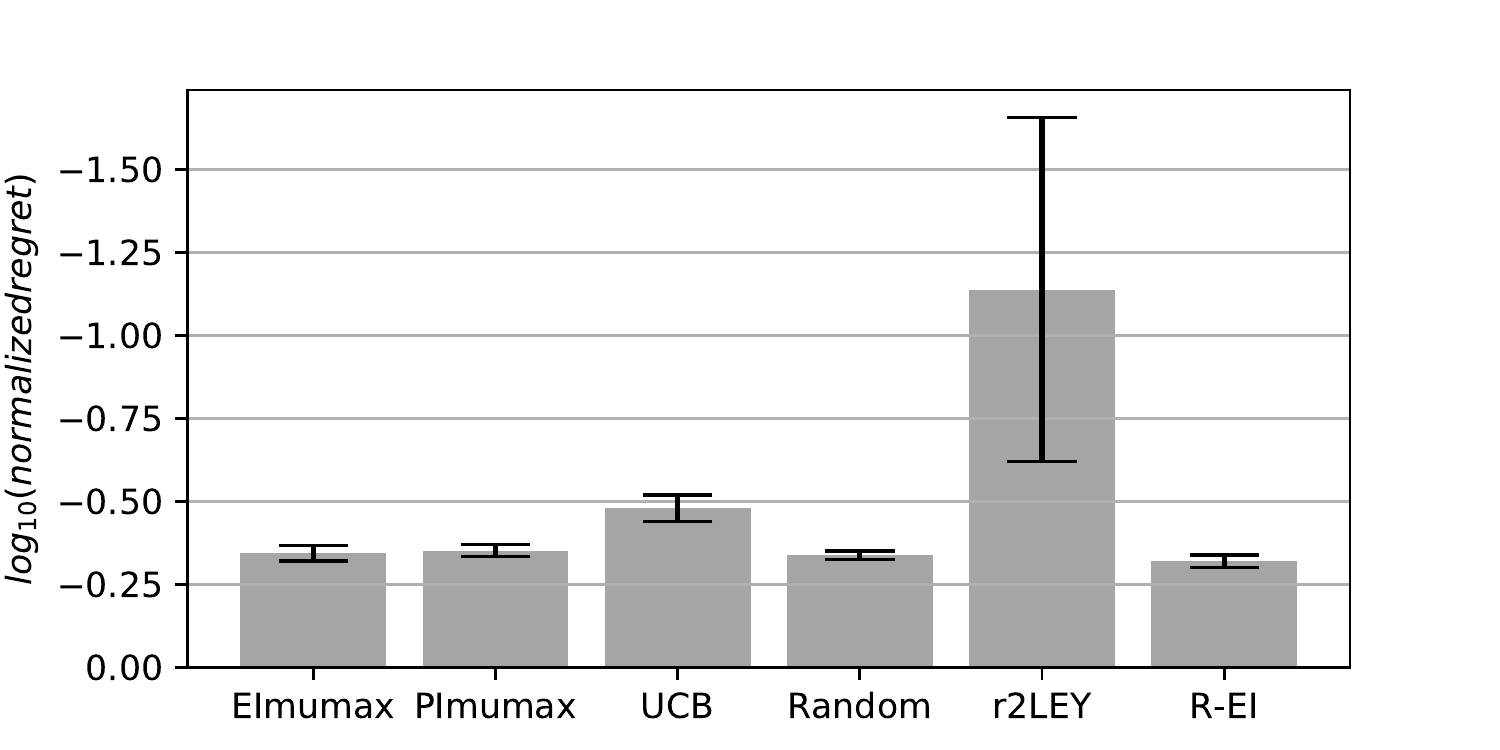}%
        \hfill
        \includegraphics[width=0.5\textwidth]{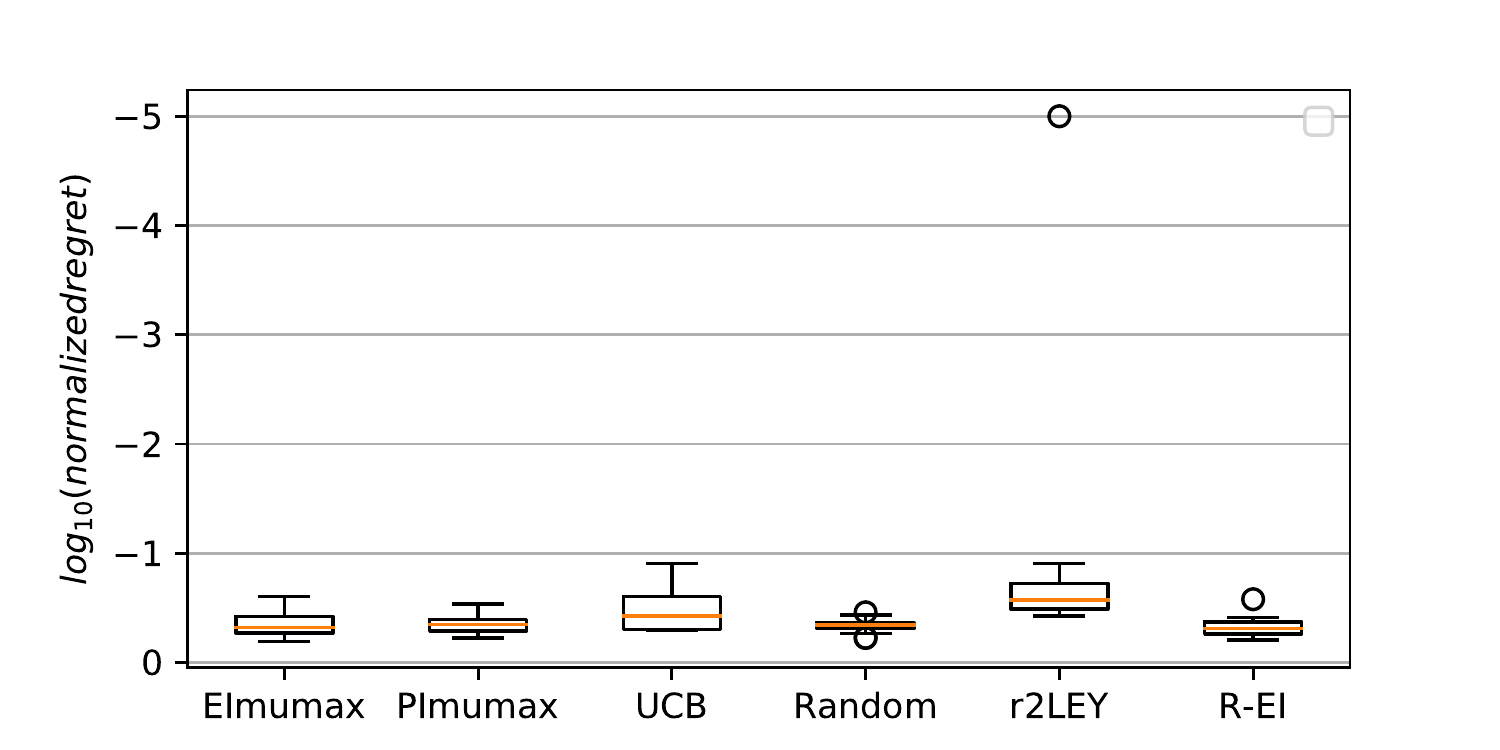}%
        \caption{Levy-8d}
        \label{sf:}
        \end{subfigure}\\ 
        \begin{subfigure}[htb!]{\textwidth}
        \includegraphics[width=0.5\textwidth]{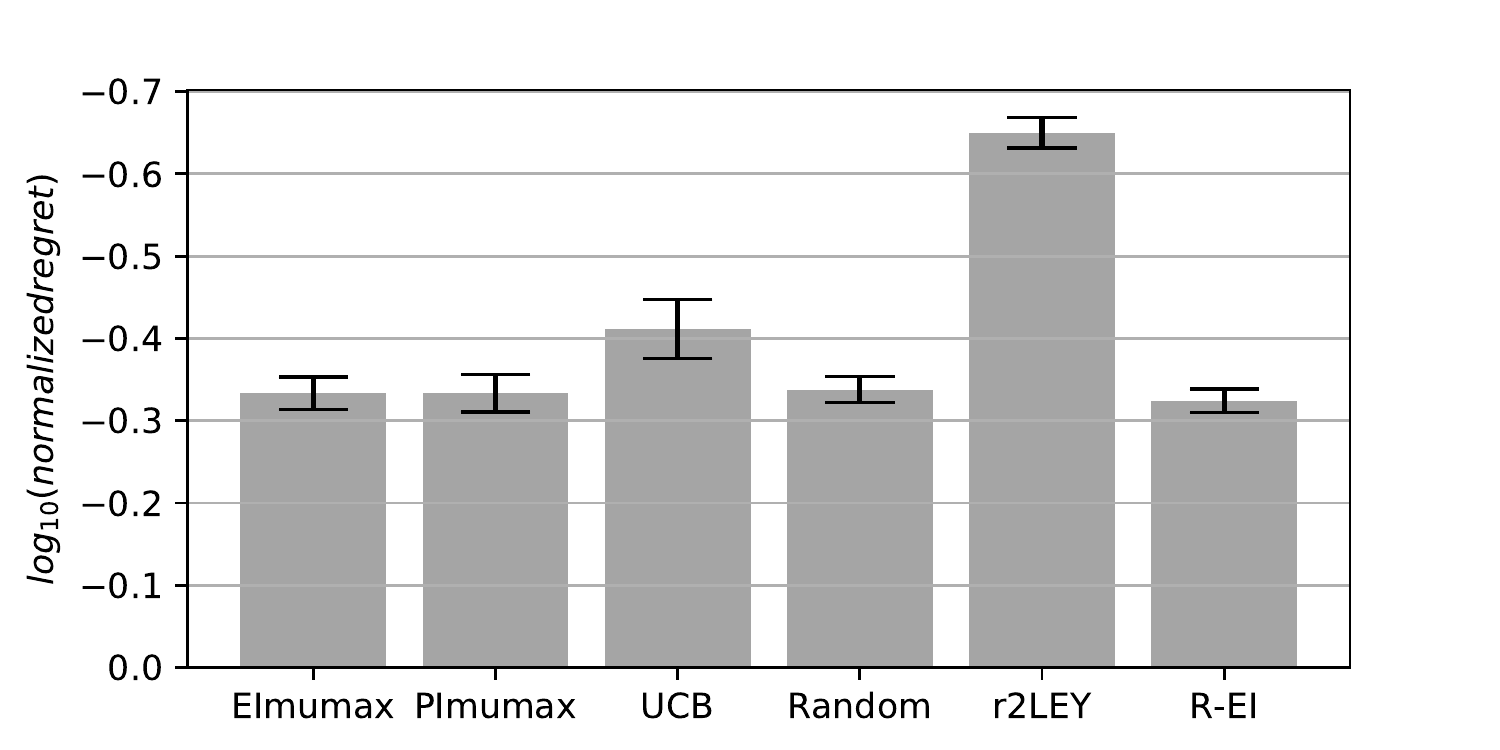}%
        \hfill
        \includegraphics[width=0.5\textwidth]{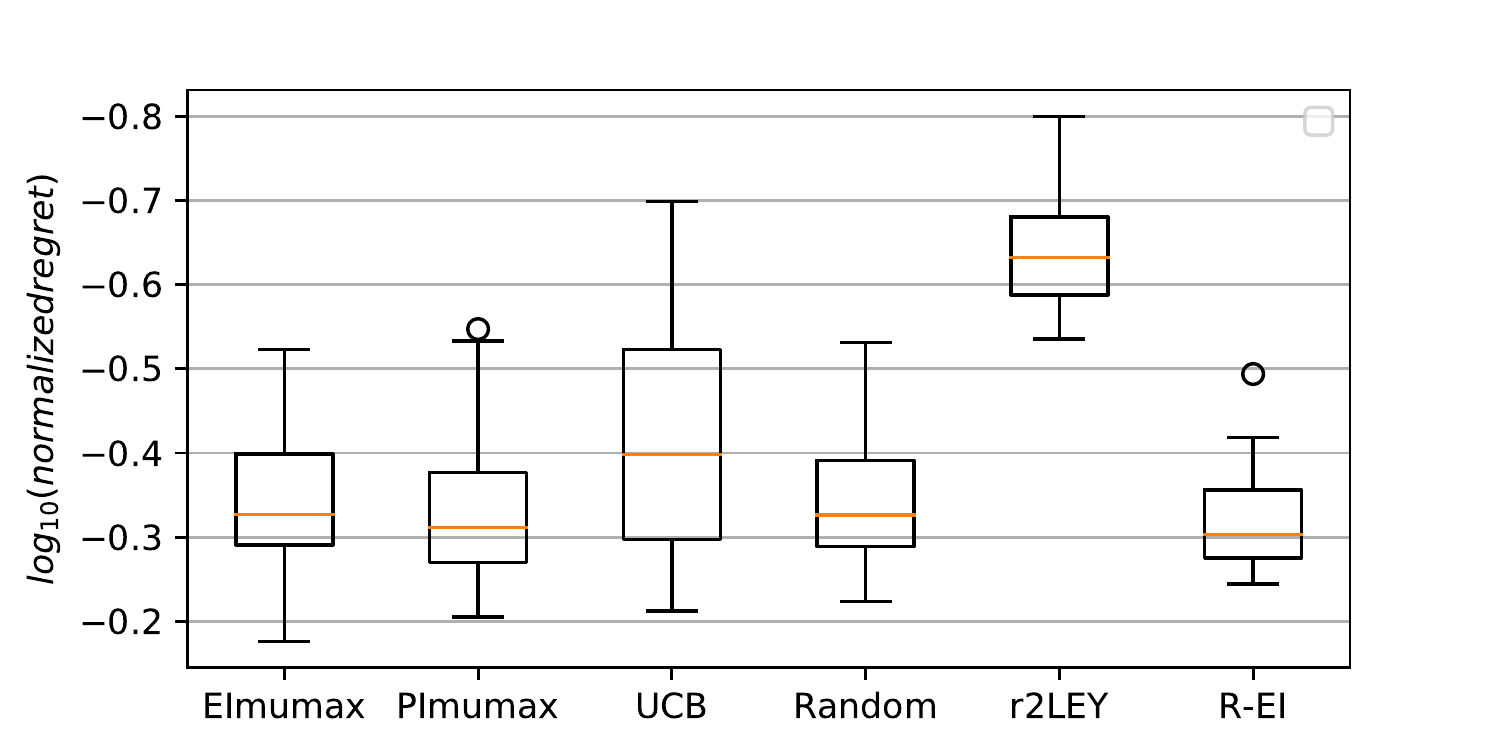}%
        \caption{Styblinski-Tang-10d}
        \label{sf:}
        \end{subfigure}\\         
    \caption{} 
    \label{f:syn_results}
\end{figure}

\begin{figure}[ht]
    \ContinuedFloat
    \centering
        \begin{subfigure}[htb!]{\textwidth}
        \includegraphics[width=0.5\textwidth]{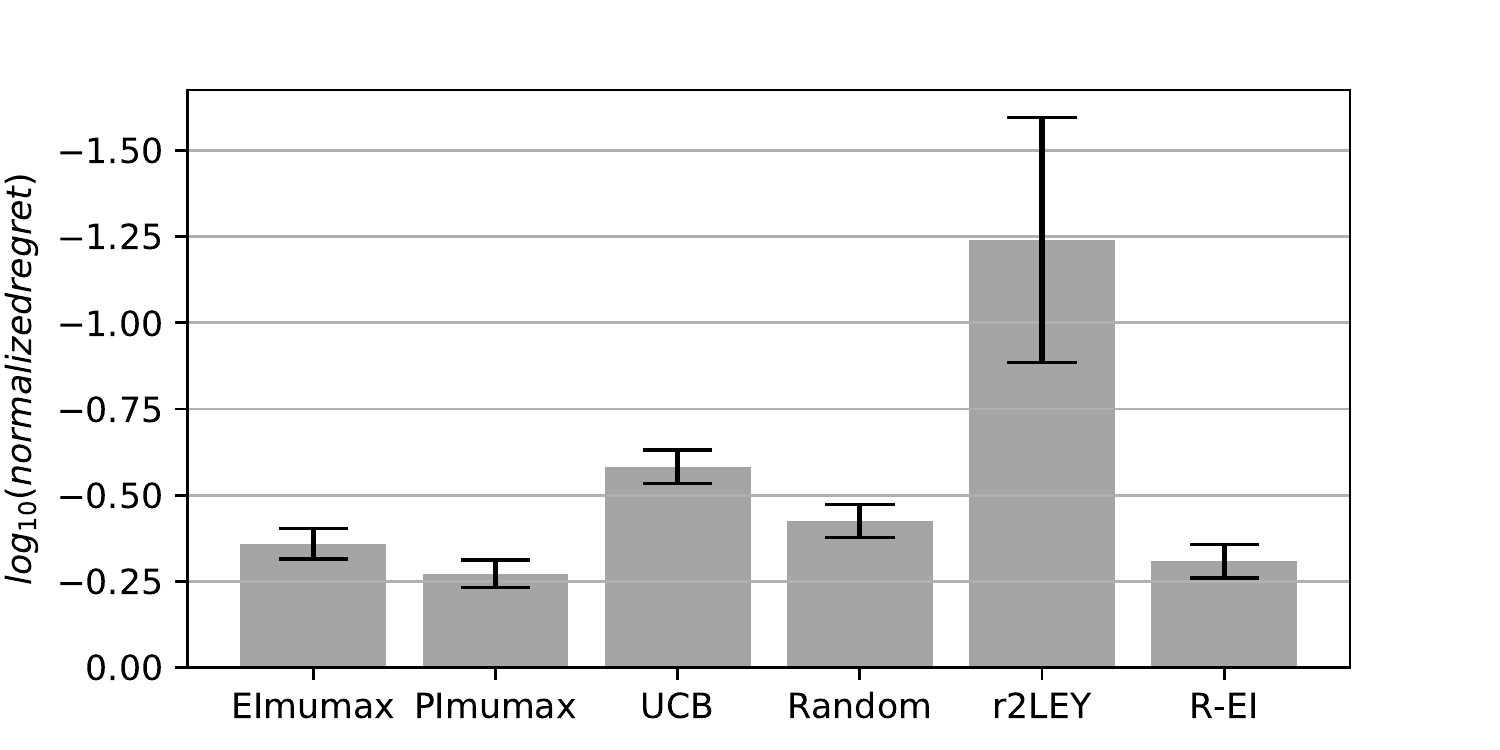}%
        \hfill
        \includegraphics[width=0.5\textwidth]{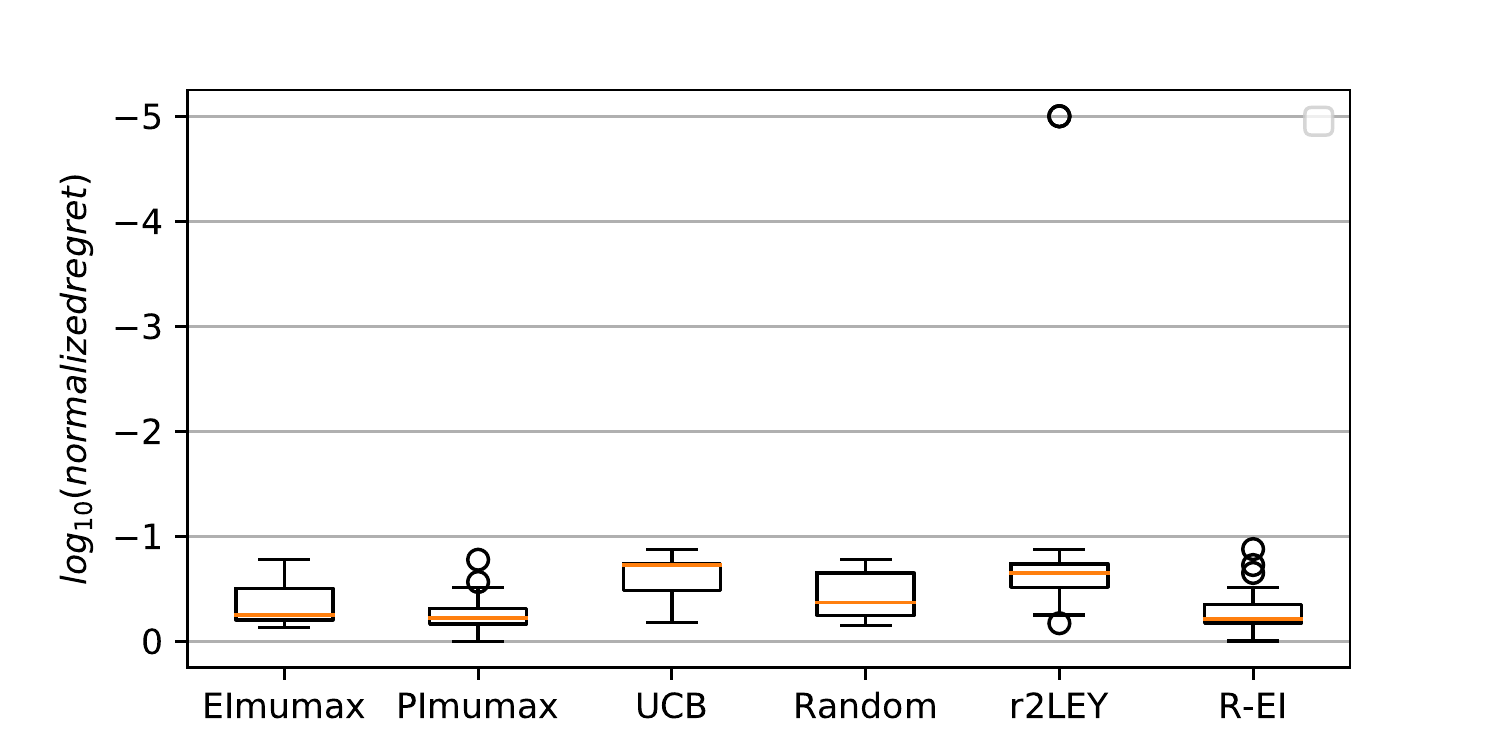}%
        \caption{Intel Sensor}
        \label{sf:}
        \end{subfigure}\\
        \begin{subfigure}[htb!]{\textwidth}
        \includegraphics[width=0.5\textwidth]{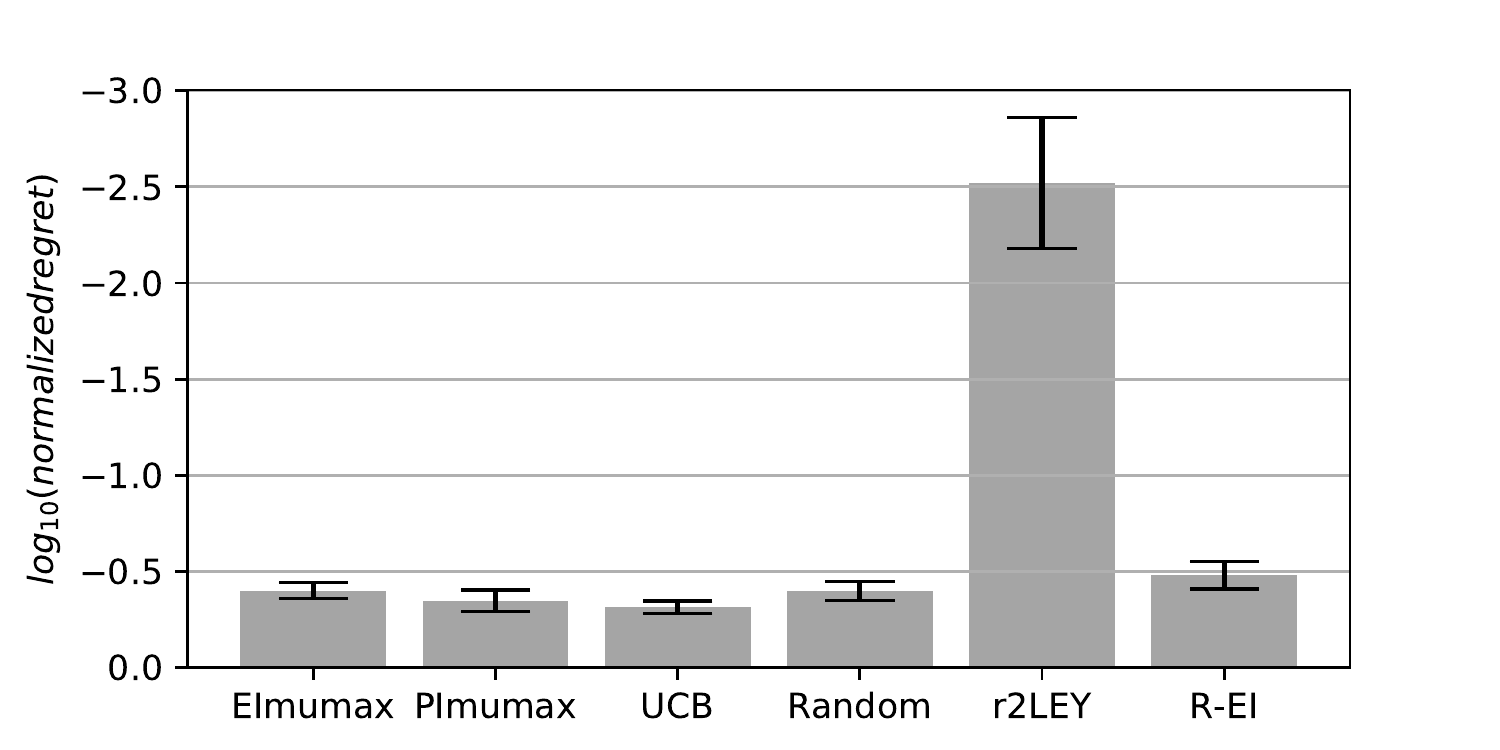}%
        \hfill
        \includegraphics[width=0.5\textwidth]{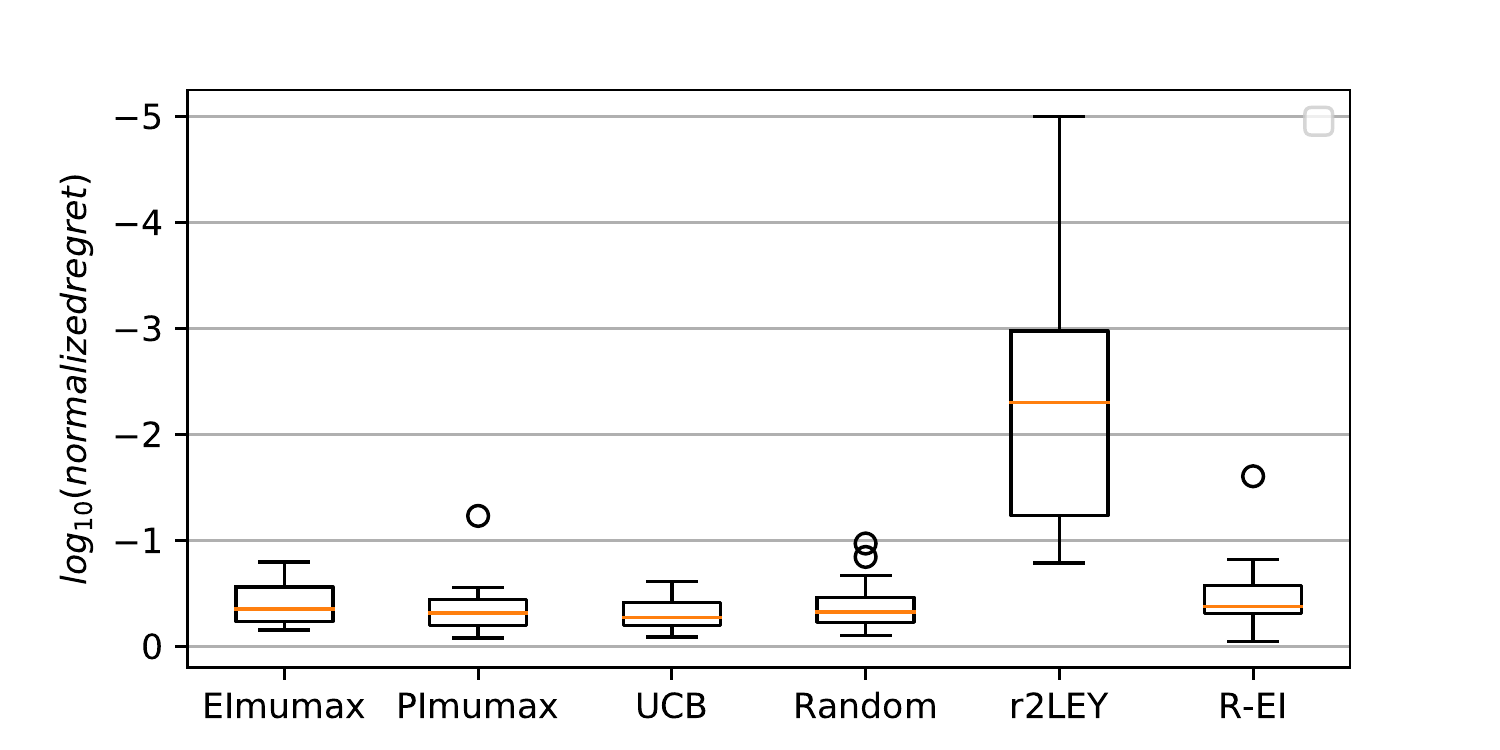}%
        \caption{SARCOS Robot}
        \label{sf:}
        \end{subfigure}\\        
    \caption{$\T{log}_{10}\T{[normalized simple regret]}$ at $T$ for all the experiments. Left column are mean $\pm$ std. error.} 
    \label{f:syn_results}
\end{figure}
\end{document}